\documentclass{amsart}
\usepackage{amssymb}

\newtheorem{theorem}{Theorem}[section]

\newtheorem{lemma}[theorem]{Lemma}
\newtheorem{corollary}[theorem]{Corollary}

\theoremstyle{definition}
\newtheorem*{definition}{Definition}

\theoremstyle{remark}

\newtheorem*{note}{Note}

\numberwithin{equation}{section}

\begin{document}
\rightline{\scriptsize CRM-2385}

\title[polynomials related to the Calogero system]
{Confluent hypergeometric orthogonal polynomials
related to the rational quantum Calogero system with harmonic confinement}

\author{J. F. van Diejen}
\address{Centre de Recherches Math\'ematiques,
Universit\'e de Montr\'eal, C.P. 6128, succursale Centre-ville,
Montr\'eal (Qu\'ebec), H3C 3J7 Canada}

\thanks{Work supported in part by
the Natural Sciences and Engineering Research Council (NSERC) of Canada.}

\subjclass{Primary 33C50; Secondary: 33C55, 33C80, 81Q05}

\date{September 1996}

\keywords{confluent hypergeometric orthogonal polynomials in
several variables, generalized spherical harmonics,
Mehta-Macdonald type integrals, integrable
quantum mechanical $n$-particle systems in dimension one}

\begin{abstract}
Two families (type $A$ and type $B$)
of confluent hypergeometric polynomials in several variables
are studied.  We describe the
orthogonality properties, differential equations, and Pieri type
recurrence formulas for these families.
In the one-variable case, the polynomials in question
reduce to the Hermite polynomials
(type $A$) and the Laguerre polynomials (type $B$), respectively.
The multivariable confluent hypergeometric families
considered here may be used to diagonalize
the rational quantum Calogero models with harmonic confinement
(for the classical root systems)
and are closely connected to the (symmetric) generalized spherical harmonics
investigated by Dunkl.
\end{abstract}

\maketitle

\section{Introduction}
In this paper multivariable orthogonal polynomials are studied associated
to the weight functions

{\em Type A (Hermite)}
\begin{subequations}
\begin{equation}\label{weightA}
\Delta^{\text{A}}(x)=\prod_{1\leq j<k\leq n} |x_j-x_k|^{2g_0}
\prod_{1\leq j\leq n} e^{-\omega x_j^2},
\end{equation}

{\em Type B (Laguerre)}
\begin{equation}\label{weightB}
\Delta^{\text{B}}(x)=\prod_{1\leq j<k\leq n} |(x_j-x_k)(x_j+x_k)|^{2g_0}
\prod_{1\leq j\leq n} |x_j|^{2g_1} e^{-\omega x_j^2} .
\end{equation}
\end{subequations}
In the one-variable case ($n=1$), the type $A$ polynomials become Hermite
polynomials ($\Delta^{\text{A}}(x)= \exp (-\omega x^2)$)
and the type $B$ polynomials reduce to Laguerre polynomials
of a quadratic argument ($\Delta^{\text{B}}(x)= |x|^{2g_1} \exp (-\omega
x^2)$).

For both families we will exhibit (systems of)
differential equations and Pieri type recurrence relations
as well as
the normalization constants that convert the polynomials
into an orthonormal system.
Multiplication of the polynomials by the square root of the weight function,
yields an orthogonal basis of eigenfunctions for the rational
quantum Calogero model with harmonic confinement
and its generalizations associated to
(the classical) root systems \cite{cal:solution,ols-per:quantum}.
The connection between this orthogonal basis and the conventional
(non-orthogonal) basis of eigenfunctions for the confined rational
Calogero model
(found by separating the quantum eigenvalue problem in
a `radial' and a `spherical' part)
is described using the theory of Dunkl's generalized
spherical harmonics with reflection group symmetry
\cite{dun:orthogonal,dun:reflection}.

The multivariable Hermite and Laguerre families
associated to the weight functions $\Delta^{\text{A}}(x)$ \eqref{weightA}
and $\Delta^{\text{B}}(x)$ \eqref{weightB}
were introduced by Macdonald \cite{mac:hypergeometric} and
Lassalle \cite{las:laguerre,las:hermite} as a generalization (more
accurately a deformation)
of the previously known special case in which
the parameter $g_0$ is being fixed at the value $1/2$
\cite{her:bessel,con:distribution,jam:special,mui:aspects}.
Recently, further insight regarding the properties of the
polynomials considered by Macdonald and Lassalle
was obtained in the context of a renewed study of the eigenvalue problem for
the rational quantum Calogero model with harmonic confinement
\cite{die:multivariable,uji-wad:rodrigues,bak-for:calogero-sutherland}.
As it turns out, some of the results reported in the present work may also be
obtained
by combining results from previous literature. For example, our evaluation
formulas
for the (squared) norms of the polynomials (cf. Theorem~\ref{normthm})
can also be gleaned from
\cite{mac:hypergeometric}, \cite{las:laguerre,las:hermite} and
\cite{bak-for:calogero-sutherland},
where expressions for these norms in a modified guise were
obtained by different methods.  (Specifically, if we make
the norm formulas in \cite{mac:hypergeometric,las:laguerre,las:hermite}
and \cite{bak-for:calogero-sutherland}
explicit with the aid of
known evaluation formulas for the Jack symmetric
functions at the identity due to Stanley \cite{sta:some,mac:symmetric}, then
they are seen
to be in correspondence with the expressions derived below in
a completely different manner.)
In all instances where overlap of this kind occurs (see the notes in
Section~\ref{sec3.4}),
our approach provides an alternative, independent
method of proof for the statements of interest.

The paper is organized as follows.
First, the confluent hypergeometric families associated to the weight functions
$\Delta^{\text{A}}(x)$ \eqref{weightA} and $\Delta^{\text{B}}(x)$
\eqref{weightB}
are defined in Section~\ref{sec2} and their main properties
(orthogonality relations, orthonormalization constants, differential equations,
and
Pieri type recurrence relations) are formulated.
In Section~\ref{sec3}, we comment in some detail
on the precise relation between our results and those obtained in previous
literature.
We will---in particular---take the opportunity to detail the
connection between the multivariable Hermite/Laguerre families
and the Calogero eigenfunctions as well as the relation to
Dunkl's generalized spherical harmonics.
In Section~\ref{sec4} we
provide the proofs for the statements in Section~\ref{sec2}
by viewing the multivariable confluent hypergeometric
Hermite  and Laguerre families of the
present work as a degeneration (viz. a limiting case)
of certain families of multivariable hypergeometric
orthogonal polynomials that were introduced in \cite{die:multivariable}
and investigated in more detail in \cite{die:properties}.
The multivariable hypergeometric polynomials relevant to us here are
generalizations of the one-variable continuous Hahn polynomials
\cite{ask:continuous,ata-sus:hahn,koe-swa:askey-scheme}
(this corresponds to type $A$)
and of the one-variable Wilson polynomials \cite{wil:some,koe-swa:askey-scheme}
(this corresponds to type $B$) to the case of several variables.
{}From a physical viewpoint,
the multivariable hypergeometric polynomials in question
are connected to the eigenfunctions of a difference (or `relativistic')
counterpart of
the rational Calogero models with harmonic confinement
\cite{die:difference,die:multivariable,rui:complete}.
The transition from
the hypergeometric to the confluent hypergeometric level corresponds
to sending the difference step-size to zero.
In this (`nonrelativistic') limit the difference Calogero model reduces to the
ordinary
Calogero model.
Some technicalities needed to perform this transition at the
level of the polynomials
(which is established by controlling the convergence of the
respective weight functions) are relegated
to an appendix at the end of the paper.

\begin{note}
Below we will always assume (unless explicitly stated
otherwise) that the parameters $g_0$ and $g_1$ entering through
the weight functions $\Delta^{\text{A}}(x)$ \eqref{weightA} and
$\Delta^{\text{B}}(x)$ \eqref{weightB} are nonnegative
and, similarly, that the scale factor $\omega$ is positive.
In principle it is possible to rescale the variables $x_1,\ldots ,x_n$ so as to
reduce
to the case that $\omega$ is fixed at the value $1$ (say). However, we have
found it
useful to keep the dependence on $\omega$ explicit in order
to have a check on the scaling properties of our expressions and so as to
suppress the
emergence of numerical constants.
\end{note}

\section{Multivariable Hermite and Laguerre polynomials}\label{sec2}
In this section the multivariable confluent hypergeometric families
associated to the weight functions $\Delta^{\text{A}}(x)$ \eqref{weightA}
and $\Delta^{\text{B}}(x)$ \eqref{weightB} are defined and the main
properties of the polynomials are stated.
The proof of these properties can be found in Section~\ref{sec4}.

\subsection{Definition and orthogonality properties}\label{sec2.1}
Let $m_\lambda$, $\lambda\in\Lambda$ denote the basis of symmetric monomials
\begin{equation}\label{monomials}
m_\lambda (x) = \sum_{\mu\in S_n (\lambda )} x_1^{\mu_1}\cdots
x_n^{\mu_n},\;\;\;\;\;\;\;
\lambda\in\Lambda
\end{equation}
with
\begin{equation}\label{cone}
\Lambda = \{ \lambda \in \mathbb{Z}^n\; |\; \lambda_1\geq\lambda_2\geq \cdots
\geq \lambda_n\geq 0 \} .
\end{equation}
Here the summation in \eqref{monomials} is over the orbit of $\lambda$ with
respect
to the action of the permutation group $S_n$ (which permutes the vector
components
$\lambda_1,\ldots ,\lambda_n$).
We will also use the notation
\begin{equation}\label{evenmonomials}
m_{2 \lambda} (x) = \sum_{\mu\in S_n (\lambda )} x_1^{2\mu_1}\cdots
x_n^{2\mu_n},\;\;\;\;\;\;\;
\lambda\in\Lambda
\end{equation}
to indicate the basis of the symmetric monomials that are even in the variables
$x_1,\ldots ,x_n$.
The monomial bases \eqref{monomials} and \eqref{evenmonomials}
inherit a partial ordering from the
dominance type partial ordering of the cone $\Lambda$ \eqref{cone} that is
defined for
$\lambda ,\mu\in\Lambda$ by
\begin{equation}\label{order}
\lambda \leq \mu \;\;\; \text{iff}\;\;\;
\sum_{1\leq j\leq k} \lambda_j \leq \sum_{1\leq j\leq k} \mu_j \;\; \text{for}
\;\; k=1,\ldots ,n
\end{equation}
($\lambda < \mu$ iff $\lambda \leq \mu$ and $\lambda \neq \mu$).

Let $\langle \cdot , \cdot \rangle_{\text{A}}$
and $\langle \cdot , \cdot \rangle_{\text{B}}$ denote the $L^2$ inner products
over $\mathbb{R}^n$
with weight function $\Delta^{\text{A}}(x)$ \eqref{weightA}
and $\Delta^{\text{B}}(x)$ \eqref{weightB}, respectively.
So, explicitly we have
\begin{subequations}
\begin{equation}\label{ipA}
\langle f , g \rangle_{\text{A}}\equiv
\int_{-\infty}^{\infty} \cdots \int_{-\infty}^{\infty}
f(x)\: \overline{g(x)}\: \Delta^{\text{A}}(x)\: dx_1\cdots dx_n
\end{equation}
(for $f,g$ in $L^2(\mathbb{R}^n,\Delta^{\text{A}}dx_1\cdots dx_n)$) and
\begin{equation}\label{ipB}
\langle f , g \rangle_{\text{B}} \equiv
\int_{-\infty}^{\infty} \cdots \int_{-\infty}^{\infty}
f(x)\: \overline{g(x)}\: \Delta^{\text{B}}(x)\: dx_1\cdots dx_n
\end{equation}
(for $f,g$ in $L^2(\mathbb{R}^n,\Delta^{\text{B}}dx_1\cdots dx_n)$).
\end{subequations}
After these notational preliminaries, we are now in the
position to define the multivariable confluent hypergeometric families
associated to the weight functions
$\Delta^{\text{A}}(x)$ \eqref{weightA}
and $\Delta^{\text{B}}(x)$ \eqref{weightB}.
\begin{definition}
The type $A$ (or multivariable Hermite) polynomials $p_\lambda^{\text{A}}(x)$,
$\lambda\in\Lambda$ are the polynomials determined (uniquely) by the conditions
\begin{itemize}
\item[A.1] $\displaystyle p_\lambda^{\text{A}}(x) = m_\lambda (x) +
\sum_{\mu\in\Lambda ,\mu < \lambda} c^{\text{A}}_{\lambda ,\mu} m_\mu (x)$,
\ \ \ \ $c^{\text{A}}_{\lambda ,\mu} \in \mathbb{C}$
\item[A.2] $\langle p_\lambda^{\text{A}} ,  m_\mu \rangle_{\text{A}} = 0$
\ if \  $\mu < \lambda$.
\end{itemize}
Similarly, the type $B$ (or multivariable Laguerre) polynomials
$p_\lambda^{\text{B}}(x)$,
$\lambda\in\Lambda$ are the polynomials determined (uniquely) by the conditions
\begin{itemize}
\item[B.1] $\displaystyle p_\lambda^{\text{B}}(x) = m_{2\lambda} (x) +
\sum_{\mu\in\Lambda ,\mu < \lambda} c^{\text{B}}_{\lambda ,\mu} m_{2\mu} (x)$,
\ \ \ \ $c^{\text{B}}_{\lambda ,\mu} \in \mathbb{C}$
\item[B.2] $\langle p_\lambda^{\text{B}} ,  m_{2\mu} \rangle_{\text{B}} = 0$
\ if \  $\mu < \lambda$.
\end{itemize}
\end{definition}

The type $A$ polynomials $p_\lambda^{\text{A}}(x)$, $\lambda\in\Lambda$
constitute a basis for the space of permutation invariant
polynomials in the variables $x_1,\ldots ,x_n$ and the type $B$ polynomials
$p_\lambda^{\text{B}}(x)$, $\lambda\in\Lambda$ form a basis for the even
subsector of this space (i.e., the subspace of symmetric polynomials
in $x_1^2,\ldots ,x_n^2$).
The following theorem says that the bases in question are orthogonal with
respect to
the inner products $\langle \cdot , \cdot \rangle_{\text{A}}$ \eqref{ipA}
and $\langle \cdot , \cdot \rangle_{\text{B}}$ \eqref{ipB}, respectively.

\begin{theorem}[Orthogonality]\label{orthothm}
Let $\lambda,\mu\in\Lambda$ \eqref{cone}. We have
\begin{eqnarray*}
\langle p_\lambda^{\text{C}} ,  p_{\mu}^{\text{C}} \rangle_{\text{C}}  =
\int_{-\infty}^{\infty} \cdots \int_{-\infty}^{\infty}
p_\lambda^{\text{C}}(x)\: \overline{p_{\mu}^{\text{C}} (x)}\:
\Delta^{\text{C}}(x)\: dx_1\cdots dx_n && \\
= 0 \;\;\; \text{if} \;\;\; \lambda \neq \mu, &&
\end{eqnarray*}
where $C$ stands for $A$ or $B$.
\end{theorem}
For weight vectors
$\lambda$ and $\mu$ that are comparable with respect to the partial order
\eqref{order}
the orthogonality of $p^{\text{C}}_\lambda$ and $p^{\text{C}}_\mu$
follows immediately from the definition of the polynomials.
Theorem~\ref{orthothm} states that the orthogonality relations in fact hold
for general weight vectors $\lambda , \mu\in \Lambda$ \eqref{cone}
(not necessarily comparable with respect to the partial
order \eqref{order}).

In order to orthonormalize the bases $\{
p_\lambda^{\text{A}}\}_{\lambda\in\Lambda}$
and $\{ p_\lambda^{\text{B}}\}_{\lambda\in\Lambda}$, it is needed to evaluate
the integrals for the (squared) norms of the polynomials.

\begin{theorem}[Norm formulas]\label{normthm}
Let $\lambda \in\Lambda$ \eqref{cone}. We have
\begin{subequations}
\begin{eqnarray*}
\langle p_\lambda^{\text{A}} ,  p_{\lambda}^{\text{A}} \rangle_{\text{A}}
\!\!\! &=& \!\!\!
\int_{-\infty}^{\infty} \cdots \int_{-\infty}^{\infty}
|p_\lambda^{\text{A}}(x)|^2\: \Delta^{\text{A}}(x)\: dx_1\cdots dx_n   \\
\!\!\! &=&\!\!\! \frac{(2 \pi)^{n/2} n!}{(2\omega )^{|\lambda
|+g_0n(n-1)/2+n/2}} \nonumber \\
\!\!\! && \!\!\! \times  \prod_{1\leq j < k\leq n}\!\!\!\!
    \frac{\Gamma ( (k-j+1)g_0 +\lambda_j-\lambda_k)\:
          \Gamma ( 1+ (k-j-1)g_0 +\lambda_j-\lambda_k)}
         {\Gamma ( (k-j)g_0 +\lambda_j-\lambda_k)\:
          \Gamma ( 1+ (k-j)g_0 +\lambda_j-\lambda_k)} \nonumber \\
\!\!\! && \!\!\! \times  \prod_{1\leq j\leq n} \Gamma (1 +(n-j) g_0 +\lambda_j
)
\end{eqnarray*}
and
\begin{eqnarray*}
\langle p_\lambda^{\text{B}} ,  p_{\lambda}^{\text{B}} \rangle_{\text{B}}
\!\!\! &=&\!\!\!
\int_{-\infty}^{\infty} \cdots \int_{-\infty}^{\infty}
|p_\lambda^{\text{B}}(x)|^2\: \Delta^{\text{B}}(x)\: dx_1\cdots dx_n  \\
\!\!\! &=& \!\!\! \frac{n!}{\omega^{2 |\lambda |+g_0n(n-1) +(g_1+1/2)n}}
\nonumber \\
\!\!\! && \!\!\! \times  \prod_{1\leq j < k\leq n}\!\!\!\!
    \frac{\Gamma ( (k-j+1)g_0 +\lambda_j-\lambda_k)\:
          \Gamma ( 1+ (k-j-1)g_0 +\lambda_j-\lambda_k)}
         {\Gamma ( (k-j)g_0 +\lambda_j-\lambda_k)\:
          \Gamma ( 1+ (k-j)g_0 +\lambda_j-\lambda_k)} \nonumber \\
\!\!\! && \!\!\! \times    \prod_{1\leq j\leq n} \Gamma (1 +(n-j) g_0
+\lambda_j )\:
                          \Gamma ( (n-j)g_0 +g_1 +1/2 +\lambda_j) \nonumber
\nonumber
\end{eqnarray*}
\end{subequations}
(where $|\lambda |\equiv \lambda_1+\cdots +\lambda_n$ and $\Gamma (\cdot )$
denotes
the gamma function).
\end{theorem}

For $\lambda=\mathbf{0}$ the
polynomials $p_\lambda^{\text{A}}$ and $p_\lambda^{\text{B}}$
reduce to the unit polynomial
($p_{\mathbf{0}}^{\text{A}}(x)=p_{\mathbf{0}}^{\text{B}}(x)=1$).
The formulas in Theorem~\ref{normthm} are then seen
to simplify to
\begin{subequations}
\begin{eqnarray}\label{mehtaA}
\langle 1 ,1\rangle_{\text{A}}
\!\!\! &=& \!\!\!
\int_{-\infty}^{\infty} \cdots \int_{-\infty}^{\infty}
\Delta^{\text{A}}(x)\: dx_1\cdots dx_n   \\
\!\!\! &=&\!\!\!
\frac{(2 \pi)^{n/2}}{(2\omega )^{g_0n(n-1)/2+n/2}}
\prod_{1\leq j\leq n} \frac{\Gamma ( 1+ jg_0)}{\Gamma (1+g_0)} \nonumber
\end{eqnarray}
and
\begin{eqnarray}\label{mehtaB}
\makebox[1em]{}\langle 1 ,1\rangle_{\text{B}}
\!\!\!\! &=& \!\!\!\!
\int_{-\infty}^{\infty} \cdots \int_{-\infty}^{\infty}
\Delta^{\text{B}}(x)\: dx_1\cdots dx_n   \\
\!\!\!\! &=&\!\!\!\!
\frac{1}{\omega^{g_0n(n-1) +(g_1+1/2)n}}\!
\prod_{1\leq j\leq n}\!\! \frac{\Gamma ( (j-1)g_0 + g_1+1/2)\, \Gamma ( 1+
jg_0)}{\Gamma (1+g_0)} .
\nonumber
\end{eqnarray}
\end{subequations}
The integrals in \eqref{mehtaA} and \eqref{mehtaB}
amount to integrals evaluated by Mehta and Macdonald
\cite{meh:random,mac:some}.
More precisely, Mehta conjectured the closed expression for
the value of the integral $\langle 1 , 1\rangle_{\text{A}}$ and Macdonald
generalized this conjecture (in terms of root systems) therewith also
including the integral $\langle 1, 1 \rangle_{\text{B}}$.
Both the integration formulas for
$\langle 1 , 1\rangle_{\text{A}}$ and $\langle 1, 1 \rangle_{\text{B}}$
in \eqref{mehtaA} and \eqref{mehtaB}
were then proven in \cite{mac:some} by viewing them
as limiting cases of an integration formula due to Selberg
\cite{sel:bemerkninger} (cf. also the introduction of
\cite{ask:some}).

\subsection{Differential equations}\label{sec2.2}
First we introduce two families of commuting difference operators
$D_{1,\beta}^{\text{A}},\ldots ,D_{n,\beta}^{\text{A}}$ and
$D_{1,\beta}^{\text{B}},\ldots ,D_{n,\beta}^{\text{B}}$.
The type $A$ operators are given by
\begin{subequations}
\begin{equation}\label{diffAbetar}
D^{\text{A}}_{r,\beta}=\!\!\!\!\!\!\!\!\!
\sum_{\stackrel{J_+,J_-\subset \{ 1,\ldots ,n\}}
               {J_+\cap J_- =\emptyset ,\; |J_+|+|J_-|\leq r}}
\!\!\!\!\!\!\!\!\! U^{\text{A}}_{J_+^c\cap J_-^c,\, r-|J_+|-|J_-|}\,
V^{\text{A}}_{J_+,J_-;\, J_+^c\cap J_-^c}\;
e^{\frac{\beta}{i}(\partial_{J_+}-\partial_{J_-})}
\end{equation}
$r=1,\ldots , n$, with
\begin{eqnarray*}
e^{\frac{\beta}{i}(\partial_{J_+}-\partial_{J_-})}\!\!\! &=&\!\!\!
\prod_{j\in J_+} e^{(\frac{\beta}{i}\frac{\partial}{\partial x_j})}
\prod_{j\in J_-} e^{-(\frac{\beta}{i}\frac{\partial}{\partial x_j})}, \\
V^{\text{A}}_{J_+, J_- ;\, K}\!\!\!  &=&\!\!\!
\prod_{j\in J_+} w^{\text{A}}(x_j)
\prod_{j\in J_-} w^{\text{A}}(-x_j)
\!\!\!\prod_{j\in J_+,j^\prime\in J_-}\!\!\!
v^{\text{A}}(x_j-x_{j^\prime})
v^{\text{A}}(x_j-x_{j^\prime}-i\beta  )\\
& & \times
\prod_{\stackrel{j\in J_+}{k\in K}}v^{\text{A}}(x_j-x_k)
\prod_{\stackrel{j\in J_-}{k\in K}}v^{\text{A}}(x_k-x_j),\\
U^{\text{A}}_{K,p}\!\!\! &=&
\!\!\!  (-1)^p \!\!\!\!\!\!\!\! \!\!\!
\sum_{\stackrel{L_+,L_-\subset K,\, L_+\cap L_-=\emptyset}
               {|L_+|+|L_-|=p}}\;
\Bigl(
\prod_{l\in L_+}\!\! w^{\text{A}}(x_l)\!\!
\prod_{l\in L_-}\!\! w^{\text{A}}(-x_l)\!\! \\
& & \makebox[7em]{}\times
\prod_{\stackrel{l\in L_+}{l^\prime\in L_-}}\!\!
v^{\text{A}}(x_l-x_{l^\prime})
v^{\text{A}}(x_{l^\prime}-x_l+i\beta )\\
& & \makebox[7em]{}\times \!\!\!\!
\prod_{\stackrel{l\in L_+}{k\in K\setminus L_+\cup L_-}} \!\!\!\!
v^{\text{A}}(x_l-x_k)\!\!\!\!
\prod_{\stackrel{l\in L_-}{k\in K\setminus L_+\cup L_-}}\!\!\!\!
v^{\text{A}}(x_k -x_l) \Bigr) ,
\end{eqnarray*}
and
\begin{equation*}
v^{\text{A}}(z)=  \Bigl(1 +\frac{\beta g_0}{iz}\Bigr) ,\;\;\;\;\;\;\;\;\;
w^{\text{A}}(z)= (1+i\beta \omega z) . \label{vwA}
\end{equation*}
The type $B$ operators read
\begin{equation}\label{diffBbetar}
D^{\text{B}}_{r,\beta}=
\sum_{\stackrel{J\subset \{ 1,\ldots ,n\} ,\, 0\leq |J|\leq r}
               {\varepsilon_j=\pm 1,\; j\in J}}
U^{\text{B}}_{J^c,\, r-|J|}\,  V^{\text{B}}_{\varepsilon J,\, J^c}\;
e^{\frac{\beta}{i}\partial_{\varepsilon J}}
\end{equation}
$r=1,\ldots , n$, with
\begin{eqnarray*}
e^{\frac{\beta}{i}\partial_{\varepsilon J}}\!\!\! &=&\!\!\!
\prod_{j\in J} e^{\varepsilon_j (\frac{\beta}{i}\frac{\partial}{\partial x_j})}
\;\;\;\;\;\;\;\;\;\;\;\;\;\;\;\;\;\;\; ( \varepsilon_j \in \{ -1 ,1\} ), \\
V^{\text{B}}_{\varepsilon J,\, K}\!\!\! &=&\!\!\!
\prod_{j\in J} w^{\text{B}}(\varepsilon_jx_j)
\prod_{\stackrel{j,j^\prime \in J}{j<j^\prime}}
v^{\text{B}}(\varepsilon_jx_j+\varepsilon_{j^\prime}x_{j^\prime})
v^{\text{B}}(\varepsilon_jx_j+\varepsilon_{j^\prime}x_{j^\prime}-i\beta )\\
& & \times
\prod_{\stackrel{j\in J}{k\in K}} v^{\text{B}}(\varepsilon_j x_j+x_k)
v^{\text{B}}(\varepsilon_j x_j -x_k),\\
U^{\text{B}}_{K,p}\!\!\! &=&\!\!\! \\
\!\!\! (-1)^p\!\!\!\!\!\!\!\!\!    & &\!\!\!\!\!
\sum_{\stackrel{L\subset K,\, |L|=p}
               {\varepsilon_l =\pm 1,\; l\in L }}\;
\Bigl( \prod_{l\in L} w^{\text{B}}(\varepsilon_l x_l)
\prod_{\stackrel{l,l^\prime \in L}{l<l^\prime}}
v^{\text{B}}(\varepsilon_lx_l+\varepsilon_{l^\prime}x_{l^\prime})
v^{\text{B}}(-\varepsilon_lx_l-\varepsilon_{l^\prime}x_{l^\prime}+i\beta )\\
& &\times
\prod_{\stackrel{l\in L}{k\in K\setminus L}} v^{\text{B}}(\varepsilon_l
x_l+x_k)
v^{\text{B}}(\varepsilon_l x_l -x_k) \Bigr) ,
\end{eqnarray*}
and
\begin{equation*}
v^{\text{B}}(z) =  \Bigl( 1 +\frac{\beta g_0}{iz} \Bigr) ,\;\;\;\;\;\;\;\;
w^{\text{B}}(z) =
\Bigl( 1 +\frac{\beta g_1}{iz}\Bigr) (1+i\beta\omega z) .
\end{equation*}
\end{subequations}
In \eqref{diffAbetar} the sum is meant over all disjoint pairs of index sets
$J_+,J_-\in \{ 1,\ldots ,n\}$ with the sum of the cardinalities
being $\leq r$. In \eqref{diffBbetar} the sum is over all index sets
$J\subset \{ 1,\ldots ,n\}$ of cardinality $\leq r$ and all
configurations of signs $\varepsilon_j\in \{ +1 ,-1\}$, $j\in J$.
We have furthermore assumed the conventions that empty products are equal to
one and
that $U_{K,p}^{\text{A}}$ and $U_{K,p}^{\text{B}}$ are equal to one for $p=0$.
The exponentials $\exp (\pm \frac{\beta}{i}\frac{\partial}{\partial x_j})$
act on analytic functions of $x_1,\ldots ,x_n$ by means of a shift of the $j$th
argument in the (for $\beta$ real) imaginary direction
\begin{equation*}
(e^{\pm (\frac{\beta}{i}\frac{\partial}{\partial x_j})}f)
(x_1,\ldots ,x_n) =
f(x_1,\ldots , x_{j-1},x_j\mp i\beta ,x_{j+1},\ldots ,x_n) .
\end{equation*}
Thus, the operators $D_{r,\beta }^{\text{A}}$ \eqref{diffAbetar}
and $D_{r,\beta }^{\text{B}}$ \eqref{diffBbetar}
are analytic difference operators of order $2r$ in
$\exp (\frac{\beta}{i}\frac{\partial}{\partial x_1}),\ldots ,
\exp (\frac{\beta}{i}\frac{\partial}{\partial x_n})$.

If we act with $D_{r,\beta }^{\text{C}}$ ($C=A,B$) on an (arbitrary) analytic
function
$f$ of the variables $x_1,\ldots ,x_n$, then we end up with an
expression that is holomorphic in the step size parameter
$\beta$ around $\beta =0$.
The coefficients of the Taylor expansion in $\beta$ around zero
can be written in terms of partial differential operators applied to the
function $f$.
We will call the first nonzero differential operator
in this expansion the leading differential part
of the difference operator. One obtains the
leading differential part $D_{r,0 }^{\text{C}}$ of $D_{r,\beta }^{\text{C}}$
($C=A,B$) by expanding the analytic difference operator
in powers of $\beta$ using the formal identity
$\exp (\pm \frac{\beta}{i}\frac{\partial}{\partial x_j})=
\sum_{m=0}^\infty (\pm\frac{\beta}{i}\frac{\partial}{\partial x_j})^m /m!$.
The result is a formal expansion of the difference operator
in terms of differential operators that
gets its precise meaning (as a Taylor expansion)
upon acting with both sides on
an (arbitrary) analytic function.

The following theorem describes the structure of the
highest order symbol of the leading
differential parts $D_{1,0 }^{\text{A}},\ldots ,D_{n,0 }^{\text{A}}$ and
$D_{1,0 }^{\text{B}},\ldots ,D_{n,0 }^{\text{B}}$
for the difference
operators $D_{r}^{\text{A}}$ \eqref{diffAbetar} and $D_{r}^{\text{B}}$
\eqref{diffBbetar},
and states that the differential operators of interest are diagonal with
respect to the
bases $\{ p_\lambda^{\text{A}}\}_{\lambda \in\Lambda}$ and
$\{ p_\lambda^{\text{B}}\}_{\lambda \in \Lambda}$, respectively.
\begin{theorem}[Differential equations]\label{diffeqthm}
{\em i.} The formal expansion of the difference operator
$D_{r,\beta }^{\text{C}}$ ($r=1,\ldots ,n$; $C=A,B$)
in the step size parameter $\beta$ has the form
\begin{equation*}
D_{r,\beta }^{\text{C}} = D_{r,0 }^{\text{C}}\:\beta^{2r} \;\;\; +
O(\beta^{2r+1})
\end{equation*}
with
\begin{equation*}
D_{r,0 }^{\text{C}} = (-1)^r\sum_{\stackrel{J\subset \{ 1,\ldots ,n\}}{|J|=r} }
\prod_{j\in J} \frac{\partial^2}{\partial x_j^2} \;\;\; + \text{l.o.}
\end{equation*}
(where l.o. stands for the terms of lower order in the derivatives
$\partial/ \partial x_1,\ldots ,\partial/ \partial x_n$).

{\em ii.} The leading differential parts $D_{1,0 }^{\text{C}}, \ldots ,D_{n,0
}^{\text{C}}$
commute and are simultaneously diagonalized by the multivariable Hermite (type
$A$) and
Laguerre (type $B$) polynomials
\begin{equation*}
D_{r,0 }^{\text{C}}\, p_\lambda^{\text{C}} = E^{\text{C}}_r(\lambda )\,
p_\lambda^{\text{C}}, \;\;\;\;\; \lambda\in\Lambda ,
\end{equation*}
with the corresponding eigenvalues given explicitly by
\begin{equation*}
E^{\text{A}}_r (\lambda )= (2\omega )^r
\sum_{\stackrel{J\subset \{ 1,\ldots ,n\}}{|J|=r}} \prod_{j\in J} \lambda_j
\;\;\;\; \text{and} \;\;\;\;
E^{\text{B}}_r (\lambda )= (4\omega )^r
\sum_{\stackrel{J\subset \{ 1,\ldots ,n\}}{|J|=r}} \prod_{j\in J} \lambda_j .
\end{equation*}
\end{theorem}
\begin{corollary}\label{limdiffcor}We have that
\begin{equation*}
\lim_{\beta\rightarrow 0} \beta^{-2r}
(D_{r,\beta }^{\text{C}}p_\lambda^{\text{C}})(x)
=E^{\text{C}}_r(\lambda )\, p_\lambda^{\text{C}}(x) , \;\;\;\;\;
\lambda\in\Lambda .
\end{equation*}
\end{corollary}
\begin{corollary}[Symmetry]\label{symmetrycor}
The leading differential operators $D_{1,0 }^{\text{C}},\ldots ,D_{n,0
}^{\text{C}}$
map the space of permutation-invariant polynomials in $x_1,\ldots ,x_n$ (type
$A$)
or $x_1^2,\ldots ,x_n^2$ (type $B$) into itself and
are symmetric with respect to the inner product
$\langle \cdot ,\cdot \rangle_{C}$ \eqref{ipA}, \eqref{ipB}
\begin{equation*}
\langle D_{r,0 }^{\text{C}} m_\lambda^{\text{C}} ,m_\mu^{\text{C}} \rangle_{C}
=
\langle m_\lambda^{\text{C}} , D_{r,0 }^{\text{C}} m_\mu^{\text{C}}
\rangle_{C},
\;\;\;\;\;\;\;\;\;\;\;\;\;\;\;\;\;\;\;\; \lambda ,\mu \in \Lambda
\end{equation*}
(with $m_\lambda^{\text{A}}(x)\equiv m_\lambda (x)$ \eqref{monomials} and
$m_\lambda^{\text{B}}\equiv m_{2\lambda}(x)$ \eqref{evenmonomials}).
\end{corollary}

The eigenvalue equations $\displaystyle (D_{r,0
}^{\text{C}}\,p_\lambda^{\text{C}})(x)=
\lim_{\beta\rightarrow 0} \beta^{-2r}
(D_{r,\beta }^{\text{C}}\,p_\lambda^{\text{C}})(x)=
E^{\text{C}}_r(\lambda )\, p_\lambda^{\text{C}}(x)$, $r=1,\ldots ,n$,
constitute a system of differential equations for the
polynomials $p_\lambda^{\text{C}}(x)$, $\lambda \in \Lambda$.
It is instructive to exhibit these differential equations in
a more explicit form
for the case $r=1$ (which amounts to the order of the equation being
equal to two).
For $r=1$, the difference operators $D_{r,\beta }^{\text{A}}$
\eqref{diffAbetar} and $D_{r,\beta }^{\text{B}}$ \eqref{diffBbetar} reduce to
\begin{subequations}
\begin{eqnarray}
\makebox[2em]{} D_{1,\beta }^{\text{A}}\!\!\!\! &=&\!\!\!\!
\sum_{1\leq j\leq n} \Bigl( w^{\text{A}}(x_j)
\prod_{1\leq k \leq n, k\neq j}\!\! v^{\text{A}}(x_j-x_k)
( e^{\frac{\beta}{i}\frac{\partial}{\partial_j}} -1 )  \\
& & \makebox[2em]{} w^{\text{A}}(-x_j)
\prod_{1\leq k \leq n, k\neq j}\!\! v^{\text{A}}(-x_j+x_k)
( e^{-\frac{\beta}{i}\frac{\partial}{\partial_j}} -1 ) \Bigr)
\nonumber
\end{eqnarray}
and
\begin{eqnarray}
\makebox[2em]{} D_{1,\beta }^{\text{B}}\!\!\! \! &=&\!\!\! \!
\sum_{1\leq j\leq n} \Bigl( w^{\text{B}}(x_j)
\prod_{1\leq k \leq n, k\neq j}\!\! v^{\text{B}}(x_j-x_k) v^{\text{B}}(x_j+x_k)
( e^{\frac{\beta}{i}\frac{\partial}{\partial_j}} -1 )  \\
& & \makebox[2em]{} w^{\text{B}}(-x_j)
\prod_{1\leq k \leq n, k\neq j}\!\! v^{\text{B}}(-x_j+x_k)
v^{\text{B}}(-x_j-x_k)
( e^{-\frac{\beta}{i}\frac{\partial}{\partial_j}} -1 ) \Bigr) \nonumber
\end{eqnarray}
\end{subequations}
(with $v^{\text{A}}(\cdot )$, $w^{\text{A}}(\cdot )$ and
$v^{\text{B}}(\cdot )$, $w^{\text{B}}(\cdot )$ taken the same as in
\eqref{diffAbetar} and \eqref{diffBbetar}, respectively).
A formal expansion in $\beta$ reveals that
$D_{1,\beta }^{\text{C}}=D_{1,0 }^{\text{C}}\: \beta^2 + O(\beta^3)$ with
\begin{subequations}
\begin{eqnarray}\label{diffA1}
\makebox[2em]{}D_{1,0}^{\text{A}}\!\!\!\!  &=&\!\!\!\!\!\!
\sum_{1\leq j\leq n} \Bigl( -\frac{\partial^2}{\partial x_j^2}
+2\omega x_j \frac{\partial}{\partial x_j} \Bigr)
-2g_0\!\!\! \sum_{1\leq j<k\leq n} \frac{1}{x_j-x_k}
\Bigl( \frac{\partial}{\partial x_j}-\frac{\partial}{\partial x_k}\Bigr) ,  \\
\makebox[2em]{}D_{1,0}^{\text{B}}\!\!\!\! &=& \!\!\!\!\!\! \label{diffB1}
\sum_{1\leq j\leq n} \Bigl( -\frac{\partial^2}{\partial x_j^2}
-2g_1 \frac{1}{x_j}\frac{\partial}{\partial x_j}
 +2\omega x_j \frac{\partial}{\partial x_j} \Bigr) \\
&&\hspace{-2em}-2g_0 \!\!\!\sum_{1\leq j<k\leq n} \left(
\frac{1}{x_j-x_k}
\Bigl( \frac{\partial}{\partial x_j}-\frac{\partial}{\partial x_k}\Bigr)
+\frac{1}{x_j+x_k}
\Bigl( \frac{\partial}{\partial x_j}+\frac{\partial}{\partial x_k} \Bigr)
\right) . \nonumber
\end{eqnarray}
\end{subequations}
Formula \eqref{diffA1} and \eqref{diffB1}
combined with the corresponding expressions for
the eigenvalues taken from Theorem~\ref{diffeqthm}
\begin{equation}\label{evAB1}
E^{\text{A}}_1 (\lambda )= 2\omega (\lambda_1+\cdots +\lambda_n), \;\;\;\;
E^{\text{B}}_1 (\lambda )= 4\omega (\lambda_1+\cdots +\lambda_n) ,
\end{equation}
render the second order differential equations
$D_{1,0}^{\text{A}}\, p^{\text{A}}_\lambda = E_1^{\text{A}}(\lambda
)p^{\text{A}}_\lambda$
and $D_{1,0}^{\text{B}}\, p^{\text{B}}_\lambda = E_1^{\text{B}}(\lambda
)p^{\text{B}}_\lambda$
in a fully explicit form.

In principle it is possible to determine for given $r$ the differential
operator $D_{r,0}^{\text{C}}$ algorithmically as the leading part of the
(explicit)
difference operator $D_{r,\beta }^{\text{C}}$
\eqref{diffAbetar} ,\eqref{diffBbetar}
(via the formal expansion in the step size parameter $\beta$).
It seems a rather nontrivial combinatorial exercise, however, to derive
along these lines a closed
expression for the differential operator $D_{r,0}^{\text{C}}$
for general $r$.

\subsection{Pieri type recurrence formulas}\label{sec2.3}
To describe the recurrence relations for the multivariable Hermite and Laguerre
polynomials
it is convenient to pass from monic polynomials to a different normalization.
Let
\begin{equation}\label{Preno}
P^{\text{A}}_\lambda (x) \equiv
c_\lambda^{\text{A}}\, p^{\text{A}}_\lambda (x),\;\;\;\;\;\;\;\;\;\;\;
P^{\text{B}}_\lambda (x) \equiv
c_\lambda^{\text{B}}\, p^{\text{B}}_\lambda (x)
\end{equation}
with
\begin{subequations}
\begin{eqnarray}\label{PrenocA}
c_\lambda^{\text{A}} &=& \prod_{1\leq j < k\leq n}
\frac{ [ (k-j)g_0 ]_{\lambda_j-\lambda_k}}
     { [ (1+k-j)g_0 ]_{\lambda_j-\lambda_k}},\\
\label{PrenocB}
c_\lambda^{\text{B}} &=& (-\omega )^{|\lambda |}
\prod_{1\leq j < k\leq n}
\frac{ [ (k-j)g_0 ]_{\lambda_j-\lambda_k}}
     { [ (1+k-j)g_0 ]_{\lambda_j-\lambda_k}} \\
&& \makebox[3em]{}\times \prod_{1\leq j\leq n}
\frac{1}{[ (n-j)g_0+g_1+1/2]_{\lambda_j}}, \nonumber
\end{eqnarray}
\end{subequations}
where we have employed the Pochhammer symbol
defined by
$[a]_0\equiv 1$ and
$[a]_l\equiv a (a+1) \cdots (a+l-1) $ for $l=1,2,3,\ldots$.

\begin{theorem}[Normalization]\label{normalizationthm}
The normalization of the polynomial $P^{\text{C}}_\lambda (x)$ \eqref{Preno} is
such that
\begin{equation*}
\lim_{\alpha\rightarrow \infty} \alpha^{-|\lambda |}
P^{\text{A}}_\lambda ( \alpha \mathbf{1}) = 1\;\;\;\;\;\;\;\;
(\text{with}\;\;\; \mathbf{1}\equiv (1,\ldots , 1))
\end{equation*}
and
\begin{equation*}
P^{\text{B}}_\lambda ( \mathbf{0} ) =1
\;\;\;\;\;\;\;\;\;\;\;\;\;\;\;\;\;\;\;\;\;\;\;\;\;
(\text{with}\;\;\; \mathbf{0}\equiv (0,\ldots , 0)).
\end{equation*}
\end{theorem}
The next theorem describes an expansion formula for the product of
$P^{\text{C}}_\lambda (x)$ \eqref{Preno} and the first elementary symmetric
function in $x_1,\ldots ,x_n$ (type $A$) or in  $x_1^2,\ldots ,x_n^2$ (type
$B$).
Formulas of this type are often referred to as Pieri formulas
\cite{mac:symmetric,sta:some}.
(More generally, Pieri formulas are relations in a commutative algebra
describing the expansion (in terms of a basis) of products
between the basis elements and a set of generators for the algebra.)
For $n=1$ the formulas in the theorem
reduce to classical three-term recurrence relations
for the one-variable Hermite and Laguerre polynomials (cf.
Section~\ref{sec3.1}).

\begin{theorem}[Pieri formulas: simplest case]\label{recrthm1}
\hfill

\noindent The (renormalized) multivariable Hermite and Laguerre polynomials
$P_\lambda^{\text{C}}(x)$ \eqref{Preno} satisfy the recurrence relations
($e_j$ denotes the $j$th unit vector in the standard basis of $\mathbb{R}^n$)
\begin{equation*}
 \bigl( \sum_{1\leq j\leq n} x_j \Bigr)  P^{\text{A}}_\lambda (x)=
\sum_{1\leq j \leq n} \Bigl(
\hat{V}_j^{\text{A}} P^{\text{A}}_{\lambda +e_j}(x) +
\hat{V}_{-j}^{\text{A}} P^{\text{A}}_{\lambda -e_j}(x) \Bigr) ,
\end{equation*}
\begin{eqnarray*}
\lefteqn{\Bigl( -\omega  \sum_{1\leq j\leq n} x_j^2 \Bigr)
    P^{\text{B}}_\lambda (x) =} && \\
&& \sum_{1\leq j \leq n} \Bigl(
\hat{V}_j^{\text{B}} P^{\text{B}}_{\lambda +e_j}(x) -
(\hat{V}_j^{\text{B}}+\hat{V}_{-j}^{\text{B}}) P^{\text{B}}_\lambda (x) +
\hat{V}_{-j}^{\text{B}} P^{\text{B}}_{\lambda -e_j}(x) \Bigr)
\end{eqnarray*}
with
\begin{eqnarray*}
\hat{V}_j^{\text{A}} &=&
\prod_{1\leq k\leq n, k\neq j}
\Bigl( 1 +\frac{g_0}{ (k-j)g_0 +\lambda_j-\lambda_k} \Bigr) ,\\
\hat{V}_{-j}^{\text{A}} &=&
\frac{ (n-j)g_0 +\lambda_j}{2\omega}
\prod_{1\leq k\leq n, k\neq j}
\Bigl( 1 -\frac{g_0}{ (k-j)g_0 +\lambda_j-\lambda_k} \Bigr)
\end{eqnarray*}
and
\begin{eqnarray*}
\hat{V}_j^{\text{B}} &=&
( (n-j)g_0 +g_1 +1/2 +\lambda_j)
\prod_{1\leq k\leq n, k\neq j}
\Bigl( 1 +\frac{g_0}{ (k-j)g_0 +\lambda_j-\lambda_k} \Bigr) ,\\
\hat{V}_{-j}^{\text{B}} &=&
((n-j)g_0 +\lambda_j)
\prod_{1\leq k\leq n, k\neq j}
\Bigl( 1 -\frac{g_0}{ (k-j)g_0 +\lambda_j-\lambda_k} \Bigr) .
\end{eqnarray*}
\end{theorem}
One word of caution is at place here. It may of course happen that for certain
$\lambda\in\Lambda$ and $j\in \{ 1,\ldots ,n\}$ the
vector $\lambda + e_j$ (or $\lambda -e_j$) does not lie in the cone $\Lambda$
\eqref{cone}.
For such boundary situations the polynomial $P_{\lambda + e_j}^{\text{C}}$
(or $P_{\lambda - e_j}^{\text{C}}$) is not defined and it might a priori seem
that r.h.s.
of the recurrence relation does not make sense in this case.
It is not difficult to verify, however, that in these situations the
coefficient
$\hat{V}_j^{\text{C}}$ (or $\hat{V}_{-j}^{\text{C}}$) in front of
$P_{\lambda + e_j}^{\text{C}}$ (or $P_{\lambda - e_j}^{\text{C}}$) vanishes.
(Indeed, for $\lambda\in\Lambda$ we have that $\lambda +e_j\not\in\Lambda$
if $\lambda_{j-1}=\lambda_j$
and that $\lambda -e_j\not\in\Lambda$ if
$\lambda_j=\lambda_{j+1}$ or if $j=n$ and $\lambda_n =0$.
In the former case we pick up a zero in $\hat{V}_j^{\text{C}}$
from the factor $1+g_0( (k-j)g_0+\lambda_j-\lambda_k)^{-1}$
with $k=j-1$ and in the latter case we have a zero in $\hat{V}_{-j}^{\text{C}}$
from the factor $1-g_0( (k-j)g_0+\lambda_j-\lambda_k)^{-1}$
with $k=j+1$ or from the factor $(n-j)g_0+\lambda_j$ with $j=n$,
respectively.)

The next step is to generalize the expansion formulas of Theorem~\ref{recrthm1}
to (multiplication by) arbitrary elementary symmetric functions.
To this end we introduce (for $r=1,\ldots ,n$)
\begin{subequations}
\begin{eqnarray}\label{esA}
\hat{E}_r^{\text{A}}(x) &=&
\sum_{\stackrel{J\subset \{ 1,\ldots ,n\}}{|J|=r}} \prod_{j\in J} x_j, \\
\label{esB}
\hat{E}_r^{\text{B}}(x) &=& (-\omega )^r
\sum_{\stackrel{J\subset \{ 1,\ldots ,n\}}{|J|=r}} \prod_{j\in J} x_j^2.
\end{eqnarray}
\end{subequations}
It is clear that the products
$\hat{E}_r^{\text{C}}(x) P_{\lambda }^{\text{C}} (x)$
can be written as a linear combination
of $P_{\mu }^{\text{C}} (x)$
with $\mu \leq \lambda +e_1+\cdots +e_r$
(this is immediate from the structure of the
monomial expansion of $P_{\lambda }^{\text{C}} (x)$ and
the fact that such expansion formulas for these products evidently hold if we
replace the
polynomials $P_{\lambda }^{\text{C}} (x)$
by their leading monomials $m_\lambda (x)$ ($C=A$) and $m_{2\lambda} (x)$
($C=B$)).
It turns out that many of the coefficients $c_\mu$ in the expansion
$\hat{E}_r^{\text{C}} P_{\lambda }^{\text{C}}= \sum_{\mu\leq \lambda} c_\mu
P_{\mu }^{\text{C}}$
are in fact zero (for $r=1$ this is of course apparent from
Theorem~\ref{recrthm1}).
The following theorem provides detailed information on the structure of
the terms entering the Pieri type expansion of the product
between the basis element $P_{\lambda }^{\text{C}} (x)$ and
an arbitrary elementary symmetric function $\hat{E}_r^{\text{C}}(x)$
\eqref{esA}, \eqref{esB}.

\begin{theorem}[Pieri formulas: general structure and leading
coefficients]\label{recrthm2}
\hfill

\noindent The (renormalized) multivariable Hermite and Laguerre polynomials
$P_\lambda^{\text{C}}(x)$ \eqref{Preno} satisfy a system
of recurrence relations of the form ($e_J \equiv \sum_{j\in J} e_j$)
\begin{eqnarray*}
\lefteqn{\hat{E}_r^{\text{C}}(x) P_\lambda^{\text{C}}(x) =} && \\
&& \sum_{\stackrel{ J_+,J_- \subset \{ 1,\ldots ,n\} ,\, \lambda
+e_{J_+}-e_{J_-}\in \Lambda}
               {J_+\cap J_- =\emptyset ,\, |J_+|+|J_-|\leq r}}
\hat{W}_{J_+ ,J_-;r}^{\text{C}}(\lambda )\:
P_{\lambda +e_{J_+}-e_{J_-}}^{\text{C}}(x) \;\;\;\;\;\;\;\;\;
r=1,\ldots ,n.
\end{eqnarray*}
The coefficients $\hat{W}_{J_+ ,J_-;r}^{\text{C}}(\lambda )$ that
correspond to index sets $J_+,J_-$ with the sum of
the cardinalities $|J_+|+|J_-|$ being equal to $r$ are explicitly given by
\begin{equation*}
\hat{W}_{J_+ ,J_-;r}^{\text{C}}(\lambda )=
\hat{V}_{J_+ ,J_-;(J_+\cup J_-)^c}^{\text{C}}
\end{equation*}
where
\begin{eqnarray*}
\hat{V}_{J_+ ,J_-;K}^{\text{A}}\!\!\!\! &=&\!\!\!\!
\prod_{j\in J_-} \Bigl( \frac{ (n-j)g_0 +\lambda_j}{2\omega} \Bigr) \\
&& \!\!\!\!\times\prod_{j\in J_+,j^\prime\in J_-}\!\!\!
\Bigl( 1 +\frac{g_0}{ (j^\prime -j)g_0 +\lambda_j-\lambda_{j^\prime}} \Bigr)
\Bigl( 1 +\frac{g_0}{ 1+(j^\prime -j)g_0 +\lambda_j-\lambda_{j^\prime}} \Bigr)
 \\
& & \!\!\!\!\times
\prod_{\stackrel{j\in J_+}{k\in K}}
\Bigl( 1 +\frac{g_0}{ (k-j)g_0 +\lambda_j-\lambda_k} \Bigr)
\prod_{\stackrel{j\in J_-}{k\in K}}
\Bigl( 1 -\frac{g_0}{ (k-j)g_0 +\lambda_j-\lambda_k} \Bigr)
\end{eqnarray*}
and
\begin{eqnarray*}
\hat{V}_{J_+ ,J_-;K}^{\text{B}}\!\!\!\! &=&\!\!\!\!
\prod_{j\in J_+} ((n-j)g_0+g_1+1/2+\lambda_j)
\prod_{j\in J_-} ((n-j)g_0+\lambda_j) \\
&&\!\!\!\! \times \prod_{j\in J_+,j^\prime\in J_-}\!\!\!
\Bigl( 1 +\frac{g_0}{ (j^\prime -j)g_0 +\lambda_j-\lambda_{j^\prime}} \Bigr)
\Bigl( 1 +\frac{g_0}{ 1+(j^\prime-j)g_0 +\lambda_j-\lambda_{j^\prime}} \Bigr)
\\
& &\!\!\!\! \times
\prod_{\stackrel{j\in J_+}{k\in K}}
\Bigl( 1 +\frac{g_0}{ (k-j)g_0 +\lambda_j-\lambda_k} \Bigr)
\prod_{\stackrel{j\in J_-}{k\in K}}
\Bigl( 1 -\frac{g_0}{ (k-j)g_0 +\lambda_j-\lambda_k} \Bigr)
\end{eqnarray*}
(with the convention that empty products are equal to one).
\end{theorem}

Theorem~\ref{recrthm2} constitutes a partial generalization of
Theorem~\ref{recrthm1}.
For $r=1$ the structure described in Theorem~\ref{recrthm2} is compatible with
that
of the formulas in Theorem~\ref{recrthm1} and we furthermore recover the
coefficient
in the r.h.s. of $P_{\lambda +e_j}^{\text{C}}$
(viz. $\hat{V}_{\{ j\} ,\emptyset ;\{ 1 ,\ldots ,n\} \setminus \{ j\}
}^{\text{C}}$)
and $P_{\lambda -e_j}^{\text{C}}$ (viz.
$\hat{V}_{\emptyset ,\{ j\} ;\{ 1 ,\ldots ,n\} \setminus \{ j\} }^{\text{C}}$)
but not that of $P_{\lambda }^{\text{C}}$
(which---according to Theorem~\ref{recrthm1}---happens to be zero for $C=A$ and
equal to
$-\sum_j (\hat{V}_{\{ j\} ,\emptyset ;\{ 1 ,\ldots ,n\} \setminus \{ j\}
}^{\text{C}} +
\hat{V}_{\emptyset ,\{ j\} ;\{ 1 ,\ldots ,n\} \setminus \{ j\} }^{\text{C}})$
for $C=B$).

Even though Theorem~\ref{recrthm2} is not completely explicit, as it does not
tell
us the expansion coefficients $\hat{W}_{J_+ ,J_-;r}^{\text{C}}(\lambda )$
for $|J_+|+|J_-|< r$, it is still useful in its present form.
For instance, the theorem implies (together with the orthogonality) that
\begin{eqnarray}\label{triple}
\lefteqn{\makebox[1em]{}
\langle \hat{E}_r^{\text{C}}\, P_\lambda^{\text{C}} , P_\mu^{\text{C}}
\rangle_{\text{C}}
=} \\
&& \left\{ \begin{array}{lll}
0 & \text{if}\;\; \mu \neq \lambda + e_{J_+}-e_{J_-} & \text{with}\;\;
|J_+|+|J_-| \leq r \\
\hat{V}_{J_+ ,J_-;(J_+\cup J_-)^c}^{\text{C}} & \text{if} \;\; \mu =\lambda +
e_{J_+}-e_{J_-} &
\text{with}\;\;  |J_+|+|J_-| = r
\end{array} \right. \nonumber
\end{eqnarray}
(where $J_+ ,J_- \in \{ 1,\ldots ,n\}$ such that $J_+ \cap J_- =\emptyset$).
When applying \eqref{triple} to the identity
$\langle \hat{E}_r^{\text{C}}\, P_\lambda^{\text{C}} , P_{\lambda
+e_{\{1,\ldots ,r\}}}^{\text{C}}
\rangle_{\text{C}}= \langle P_\lambda^{\text{C}} , \hat{E}_r^{\text{C}}\,
P_{\lambda +e_{\{1,\ldots ,r\}}}^{\text{C}}\rangle_{\text{C}}$ one arrives at a
system of recurrence
relations for the squared norm of $P_\lambda^{\text{C}}$
\begin{eqnarray}\label{recrnorm}
\lefteqn{\hat{V}_{\{ 1,\dots ,r\} ,\emptyset ;\{ r ,\ldots ,n\}
}^{\text{C}}(\lambda)
\:
 \langle P_{\lambda +e_{\{1,\ldots ,r\}}}^{\text{C}}, P_{\lambda +e_{\{1,\ldots
,r\}}}^{\text{C}}
\rangle_{\text{C}} =} && \\
&& \hat{V}_{\emptyset ,\{ 1,\dots ,r\} ;\{ r ,\ldots ,n\} }^{\text{C}}(\lambda
+e_{\{ 1,\ldots ,r\}})
\: \langle P_{\lambda }^{\text{C}}, P_{\lambda }^{\text{C}}
\rangle_{\text{C}} \nonumber
\end{eqnarray}
($r=1,\ldots ,n$). The recurrence relations in \eqref{recrnorm} determine
$\langle P_\lambda^{\text{C}}, P_\lambda^{\text{C}} \rangle_{\text{C}}$
uniquely in terms
of $\langle 1,1 \rangle_{\text{C}}$
(because the (fundamental weight)
vectors $e_{\{ 1,\ldots ,r\}}$, $r=1,\ldots ,n$ positively generate the cone
$\Lambda$ \eqref{cone} and the coefficient
$\hat{V}_{\{ 1,\dots ,r\} ,\emptyset ;\{ r ,\ldots ,n\} }^{\text{C}}(\lambda)
\neq 0$
for $\lambda\in\Lambda$).
This observation gives rise to an alternative (constructive) proof of the norm
formulas in Theorem~\ref{normthm} different from the proof presented in
Section~\ref{sec4.1}.
Indeed, by using the property
$c_\lambda^{\text{C}}=c_{\lambda +e_{\{ 1,\ldots ,r\}}}^{\text{C}}
\hat{V}_{\{ 1,\dots ,r\} ,\emptyset ;\{ r ,\ldots ,n\} }^{\text{C}}(\lambda)$
one rewrites \eqref{recrnorm} in the monic form
\begin{eqnarray}
\lefteqn{\langle p_{\lambda +e_{\{1,\ldots ,r\}}}^{\text{C}},
p_{\lambda +e_{\{1,\ldots ,r\}}}^{\text{C}}\rangle_{\text{C}}=} && \\
&& \hat{V}_{\{ 1,\dots ,r\} ,\emptyset ;\{ r ,\ldots ,n\} }^{\text{C}}(\lambda)
\hat{V}_{\emptyset ,\{ 1,\dots ,r\} ;\{ r ,\ldots ,n\} }^{\text{C}}(\lambda
+e_{\{ 1,\ldots ,r\}})
\,\langle p_\lambda^{\text{C}}, p_\lambda^{\text{C}} \rangle_{\text{C}}
,\nonumber
\end{eqnarray}
which upon iteration and matching of the initial conditions so as to reduce
for $\lambda = \mathbf{0}$ to
the Mehta-Macdonald formulas \eqref{mehtaA}, \eqref{mehtaB}
(cf. Section~\ref{sec4.1}) leads to the
norm formulas of Theorem~\ref{normthm}.

The last theorem (below) provides a complete (explicit) description
of the expansion coefficients $\hat{W}_{J_+ ,J_-;r}^{\text{C}}(\lambda )$
for the type $B$ case
(thus including the coefficients corresponding to index sets with
$|J_+|+|J_-|<r$).
This renders the system of recurrence relations of the type
given by Theorem~\ref{recrthm2}
in a fully explicit form for the multivariable Laguerre family.

\begin{theorem}[Pieri formulas: explicit expansion coefficients Laguerre
case]\label{recrthm3}
\hfill

\noindent The coefficients in the recurrence relations of the type described by
Theorem~\ref{recrthm2}
are for the renormalized multivariable Laguerre polynomials
$P_\lambda^{\text{B}}(x)$ \eqref{Preno} given by
\begin{equation*}
\hat{W}_{J_+ ,J_-;r}^{\text{B}}(\lambda )=
\hat{V}^{\text{B}}_{J_+,J_-;\, (J_+\cup J_-)^c}\,
\hat{U}^{\text{B}}_{(J_+\cup J_-)^c,\, r-|J_+|-|J_-|}
\end{equation*}
with $\hat{V}^{\text{B}}_{J_+,J_-;\, (J_+\cup J_-)^c}$ taken the same
as in Theorem~\ref{recrthm2} and
\begin{eqnarray*}
&& \hat{U}^{\text{B}}_{K,\, p}= (-1)^p \times \\
&& \!\!\!\!\!\!\!\! \!\!
\sum_{\stackrel{L_+,L_-\subset K,\, L_+\cap L_-=\emptyset}
               {|L_+|+|L_-|=p}}\;
\Biggl(
\prod_{l\in L_+} ((n-l)g_0+g_1+1/2+\lambda_j)
\prod_{l\in L_-} ((n-l)g_0+\lambda_j) \\
&&\!\!\!\! \times \prod_{l\in L_+,l^\prime\in L_-}\!\!\!
\Bigl( 1 +\frac{g_0}{ (l^\prime -l)g_0 +\lambda_l-\lambda_{l^\prime}} \Bigr)
\Bigl( 1 -\frac{g_0}{ 1+(l^\prime-l)g_0 +\lambda_l-\lambda_{l^\prime}} \Bigr)
\\
& &\!\!\!\! \times
\prod_{\stackrel{l\in L_+}{k\in K}}
\Bigl( 1 +\frac{g_0}{ (k-l)g_0 +\lambda_l-\lambda_k} \Bigr)
\prod_{\stackrel{l\in L_-}{k\in K}}
\Bigl( 1 -\frac{g_0}{ (k-l)g_0 +\lambda_l-\lambda_k} \Bigr)  \Biggr)
\end{eqnarray*}
(with the convention that $\hat{U}^{\text{B}}_{K,\, p}\equiv 1$
for $p=0$).
\end{theorem}

\section{Comments}\label{sec3}

\subsection{The special case $n=1$}\label{sec3.1}
In the case of one single variable ($n=1$), the weight functions
reduce to
\begin{equation}
\Delta^{\text{A}} (x) = e^{-\omega x^2},
\;\;\;\;\;\;\;\;\;\;\;\;\;\;\;\;\;\;\;\;\;\;
\Delta^{\text{B}} (x) = |x|^{2g_1}\,e^{-\omega x^2} .
\end{equation}
The polynomials then become monic Hermite polynomials (type $A$) and
monic Laguerre polynomials of a quadratic argument (type $B$),
which can be written explicitly in terms of a terminating
confluent hypergeometric series \cite{abr-ste:handbook}
\begin{subequations}
\begin{eqnarray}\label{exprepA}
p_\lambda^{\text{A}}(x)\!\!\! &=&\!\!\! \left\{
\begin{array}{ll}
\frac{[1/2]_{\lambda /2}}{(-\omega )^{\lambda /2}}
\; {}_1F_1 \left( \begin{array}{c} -\lambda /2 \\ 1/2 \end{array}\, ;\,
\omega x^2 \right) & \text{for}\; \lambda\; \text{even} \\
  \frac{[3/2]_{(\lambda -1) /2}}{(-\omega )^{(\lambda -1) /2}}
\; x\;
{}_1F_1 \left( \begin{array}{c} -(\lambda -1) /2 \\ 3/2 \end{array}\, ;\,
\omega x^2 \right) & \text{for}\; \lambda\; \text{odd}
\end{array} \right. \\
p_\lambda^{\text{B}}(x)\!\!\! &=& \!\!\!
\frac{[g_1 +1/2]_{\lambda}}{(-\omega )^\lambda}\;
{}_1F_1 \left( \begin{array}{c} -\lambda \\ g_1 + 1/2 \end{array}\, ;\,
\omega x^2 \right) \label{exprepB}
\end{eqnarray}
\end{subequations}
(with $\displaystyle {}_1F_1
\left( \begin{array}{c}  a\\ b \end{array} ; z \right)
\equiv \sum_{m=0}^\infty \frac{[a]_m}{[b]_m m!} z^m$ and
$\lambda =0,1,2,\ldots$).
The norm formulas, differential equations and
recurrence relations reduce in this special situation to
classical formulas for the one-variable Hermite and the Laguerre polynomials
(notice however that the scale parameter $\omega$ is usually
taken to be equal to one and that
our normalization differs from the standard one)

\noindent {\em Norm formulas} (cf. Theorem~\ref{normthm})
\begin{subequations}
\begin{eqnarray}
\int_{-\infty}^\infty | p_\lambda^{\text{A}}(x)|^2 \Delta^{\text{B}}(x)\, dx
&=&\frac{\lambda ! \sqrt{\pi}}{2^\lambda \omega^{\lambda +1/2}} \\
\int_{-\infty}^\infty | p_\lambda^{\text{B}}(x)|^2 \Delta^{\text{B}}(x)\, dx
&=&\frac{\lambda ! \Gamma (g_1 +1/2 +\lambda)}
        {\omega^{2\lambda +g_1+1/2}}
\end{eqnarray}
\end{subequations}
{\em Differential equations} (cf. Theorem~\ref{diffeqthm} and
Eqs.~\eqref{diffA1}, \eqref{diffB1} and \eqref{evAB1})
\begin{subequations}
\begin{eqnarray}
-\frac{d^2}{dx^2} p_\lambda^{\text{A}}
+2\omega x \frac{d}{dx} p_\lambda^{\text{A}} =
2\omega\lambda \:  p_\lambda^{\text{A}} \\
-\frac{d^2}{dx^2} p_\lambda^{\text{B}}
-\frac{2g_1}{x}\frac{d}{dx} p_\lambda^{\text{B}}
+2\omega x \frac{d}{dx} p_\lambda^{\text{B}} =
4\omega\lambda\:  p_\lambda^{\text{B}} \label{lagdiffB}
\end{eqnarray}
\end{subequations}
{\em Recurrence relations} (cf. Theorem~\ref{recrthm1})
\begin{subequations}
\begin{eqnarray}
x P^{\text{A}}_\lambda  &=&
 P^{\text{A}}_{\lambda +1}  +
\frac{\lambda}{2\omega}  P^{\text{A}}_{\lambda -1}  \\
-\omega x^2 P^{\text{B}}_\lambda  &=&
(g_1+1/2+\lambda ) P^{\text{B}}_{\lambda +1}
-(g_1+1/2+2\lambda ) P^{\text{B}}_{\lambda }
+\lambda P^{\text{B}}_{\lambda -1}
\end{eqnarray}
\end{subequations}
with $P_\lambda^{\text A}(x)=p^{\text A}_\lambda(x)$ and
$P_\lambda^{\text B}(x)=
\frac{(-\omega )^\lambda}{[g_1+1/2]_\lambda} p^{\text B}_\lambda(x)$.
The normalization properties for
$P_\lambda^{\text A}(x)$ and $P_\lambda^{\text B}(x)$
in Theorem~\ref{normalizationthm} state that
$\lim_{\alpha \rightarrow \infty}
\alpha^{-\lambda} P_\lambda^{\text A} (\alpha )=1$ and that
$P_\lambda^{\text B}(0)=1$.
In the present situation these properties are immediate
from the explicit confluent hypergeometric
representations in \eqref{exprepA} and \eqref{exprepB}.

\subsection{Eigenfunctions for the rational Calogero model in an harmonic
well}\label{sec3.2}
Most of the results stated in Section~\ref{sec2}
admit an interpretation in terms of certain
exactly solvable quantum mechanical $n$-particle models on the line.
Specifically, conjugation with the square root of the weight function
$\Delta^{\text{C}}(x)$ \eqref{weightA}, \eqref{weightB}
transforms the second order differential operators
$D_1^{\text{C}}$ \eqref{diffA1}, \eqref{diffB1}
into Hamiltonians for the rational quantum Calogero models
with harmonic confinement associated
to the classical root systems \cite{cal:solution,ols-per:quantum}
\begin{subequations}
\begin{eqnarray}\label{HA1}
\lefteqn{H^{\text{A}}_1 =
(\Delta^{\text{A}})^{\frac{1}{2}} D^{\text{A}}_1
(\Delta^{\text{A}})^{-\frac{1}{2}} =} && \\
&& \sum_{1\leq j\leq n}
\Bigl(- \frac{\partial^2}{\partial x_j^2} + \omega^2 x_j^2 \Bigr)
+2g_0(g_0-1) \sum_{1\leq j<k\leq n} (x_j-x_k)^{-2}
\;\;\; -E_0^{\text{A}} ,\nonumber
\end{eqnarray}
\begin{eqnarray}\label{HB1}
\lefteqn{H^{\text{B}}_1 =
(\Delta^{\text{B}})^{\frac{1}{2}} D^{\text{B}}_1
(\Delta^{\text{B}})^{-\frac{1}{2}} =} && \\
&& \sum_{1\leq j\leq n}
\Bigl( -\frac{\partial^2}{\partial x_j^2} + g_1(g_1-1)x_j^{-2}+
\omega^2 x_j^2 \Bigr)  \nonumber \\
&&  +2g_0(g_0-1) \sum_{1\leq j<k\leq n}
\Bigl( (x_j-x_k)^{-2} +(x_j-x_k)^{-2} \Bigr)
\;\;\; -E_0^{\text{B}} \nonumber
\end{eqnarray}
\end{subequations}
with $E_0^{\text{A}}=\omega n (1+g_0(n-1))$ and
$E_0^{\text{B}}=\omega n (1+2g_0(n-1)+2g_1)$.
It is clear from Theorem~\ref{diffeqthm} that the functions
\begin{equation}\label{psiAB}
\psi_\lambda^{\text{C}} (x) =
\left(\Delta^{\text{C}}(x)\right)^{\frac{1}{2}}\, p^{\text{C}}_\lambda (x) ,
\;\;\;\; \lambda\in\Lambda\;\;\;\;\;\;\;\;\;(C=A,B)
\end{equation}
constitute a basis of eigenfunctions for $H_1^{\text{C}}$ \eqref{HA1},
\eqref{HB1}
with the corresponding eigenvalues given by $E_1^{\text{C}}(\lambda )$
\eqref{evAB1}.
In its present form these eigenfunctions for the rational Calogero
model with harmonic term were introduced in
\cite{die:multivariable,bak-for:calogero-sutherland} and also
(for type $A$) in \cite{uji-wad:rodrigues}.
The orthogonality (in $L^2(\mathbb{R}^n,d x_1,\ldots ,dx_n)$)
of the basis $\psi_\lambda^{\text{C}} (x)$ \eqref{psiAB}
follows from Theorem~\ref{orthothm}
and the orthonormalization constants can
be read-off from Theorem~\ref{normthm}.

Historically, the study of the eigenvalue problem for the type $A$ Hamiltonian
\eqref{HA1} was initiated by Calogero,
who computed the spectrum and determined the structure of
the corresponding eigenfunctions to be a product of the ground-state wave
function
$\psi_{\mathbf{0}}^{\text{A}}(x)=(\Delta^{\text{A}}(x))^{1/2}$
and certain symmetric polynomials in $x_1,\ldots ,x_n$ \cite{cal:solution}.
To be precise, Calogero considered a translationally symmetric
$n$-particle system with a potential of the form
\begin{equation*}
V(x) = \sum_{1\leq j< k\leq n}  \bigl(\, G_0 (x_j-x_k)^{-2} + G_1
(x_j-x_k)^{2}\, \bigr) ,
\end{equation*}
which is seen to be equivalent to the type $A$ system above up to a simple
center of mass motion ($\sum_{j<k} (x_j-x_k)^2 = n\sum_j x_j^2-(\sum_j
x_j)^2$).
Explicit expressions for the Calogero eigenfunctions
were found for small particle number ($n\leq 5$)
by Perelomov and Gambardella \cite{per:algebraic,gam:exact}.
More recently, it was observed
\cite{bri-han-vas:explicit,bri-han-kon-vas:calogero}
that for arbitrary particle number $n$ it is possible to construct a basis of
eigenfunctions for $H_1^{\text{A}}$ \eqref{HA1}
with the aid of certain creation and annihilation operators
$a_j^+$ and $a_j^-$ ($j=1,\ldots ,n$) that are built from the
differential-reflection operators
introduced by Dunkl \cite{dun:differential-difference}.
In a nutshell: the relevant operators $a_j^\pm$ are obtained starting from
the usual creation/annihilation operators
$(\mp\frac{\partial}{\partial x_j} + \omega x_j)$ for a system
of uncoupled (bosonic) harmonic oscillators (this corresponds to $g_0=0$)
by replacing the partial derivatives
by the corresponding Dunkl differential-reflection operators (associated to the
root system $A_{n-1}$).
If one acts on the ground-state wave function
$\Psi_{\mathbf{0}}^{\text{A}}(x)= (\Delta^{\text{A}} (x))^{1/2}$
(which is annihilated by $a_1^-,\ldots ,a_n^-$)
with an arbitrary symmetric homogeneous polynomial of degree $l$
in the creation operators $a_1^+,\ldots ,a_n^+$, then one
winds up with a symmetric eigenfunction of $H^{\text{A}}$ \eqref{HA1}
with eigenvalue $2l\omega$.
By taking for the symmetric polynomial in question the Jack
symmetric function $J_\lambda (x;1/g_0)$ \cite{sta:some,mac:symmetric},
one arrives (up to a normalization factor) precisely at the eigenfunction
$\Psi_\lambda^{\text{A}}(x)$ \eqref{psiAB} \cite{uji-wad:rodrigues}.

The basis of eigenfunctions of the form $\Psi_\lambda^{\text{C}}(x)$
\eqref{psiAB}
for $H_1^{\text{C}}$ \eqref{HA1}, \eqref{HB1} is also very special in that it
simultaneously diagonalizes the higher-order quantum integrals
for the Calogero model. Specifically, by applying
the similarity transformation
$D_1^{\text{C}}\rightarrow H_1^{\text{C}}$ also to the higher-order
differential
operators $D_r^{\text{C}}$ in Theorem~\ref{diffeqthm},
one obtains a complete set of commuting quantum
integrals for the Calogero model of the form
\begin{eqnarray}\label{HABr}
\lefteqn{H^{\text{C}}_r =
(\Delta^{\text{C}})^{\frac{1}{2}} D^{\text{A}}_r
(\Delta^{\text{C}})^{-\frac{1}{2}} =} \\
&&\makebox[3em]{} (-1)^r \sum_{\stackrel{J\in \{ 1,\ldots ,n\}}{|J|=r}}
\prod_{j\in J} \frac{\partial^2}{\partial x_j^2}
\;\;\; +\text{l.o.}
\;\;\;\;\;\;\;\;\;\; (r=1,\ldots ,n). \nonumber
\end{eqnarray}
It is immediate from Theorem~\ref{diffeqthm} that the functions
$\Psi_\lambda^{\text{C}}(x)$ \eqref{psiAB} constitute
a basis of joint eigenfunctions for the differential operators
$H^{\text{C}}_r$ \eqref{HABr}, $r=1,\ldots ,n$.
Since the corresponding eigenvalues
$E_1^{\text{C}}(\lambda ),\ldots ,E_n^{\text{C}}(\lambda)$ separate
the points of the cone $\Lambda$ \eqref{cone}, it follows that
this property actually determines such basis uniquely.
In other words,
ambiguities in the possible choices for the basis of eigenfunctions
of $H_1^{\text{C}}$ \eqref{HA1}, \eqref{HB1} caused by the degeneracy
of the spectrum \eqref{evAB1} get eliminated by requiring
the eigenfunctions to be a basis of joint eigenfunctions
for $H_1^{\text{C}},\ldots ,H_n^{\text{C}}$.
The other important (and restrictive) property of the basis
$\Psi_\lambda^{\text{C}}(x)$ \eqref{psiAB},
viz. its orthogonality, implies---together with
the real-valuedness of the eigenvalues $E_r^{\text{C}}(\lambda)$---that the
differential
operators $H_r^{\text{C}}$ are (essentially) self-adjoint
(since unitarily equivalent to real multiplication operators)
in the Hilbert space of permutation-invariant (type $A$) or
permutation-invariant and even (type $B$)
functions in $L^2 (\mathbb{R}^n,dx_1\cdots dx_n)$
(cf. Corollary~\ref{symmetrycor}).

One would like to cast the differential operators
$H_r^{\text{C}}$ for arbitrary
$r\in \{ 1,\ldots ,n\}$ in a more explicit form.
Explicit formulas for $n$-independent commuting differential operators
generating the same algebra as our operators
$H_1^{\text{B}},\ldots ,H_n^{\text{B}}$ appear as special cases
of the commuting families presented
in \cite{och-osh-sek:commuting,osh-sek:commuting}.
In order to
determine the relation between the operators $H_r^{\text{B}}$ \eqref{HABr}
and the type $B$ differential operators
in \cite{och-osh-sek:commuting,osh-sek:commuting} explicitly
(for all $r$), one would have to know the eigenvalues of the latter
operators on the basis $\Psi_\lambda^{\text{B}} (x)$ \eqref{psiAB}.
For the type $A$ we are not aware that similar explicit formulas
for a set of commuting differential operators generating the same
(commuting) algebra as $H_1^{\text{A}},\ldots ,H_n^{\text{A}}$
have been reported in the literature (except for $\omega =0$
\cite{och-osh-sek:commuting}).
However, it is known that such commuting differential operators
may be characterized in terms of Dunkl's
differential-reflection operators \cite{pol:exchange}
and recently the eigenvalues of the thus obtained differential
operators on the basis $\Psi_\lambda^{\text{A}} (x)$ \eqref{psiAB}
were determined \cite{kak:common} (thus allowing a comparison
with our differential operators $H_r^{\text{A}}$).

\subsection{Relation to Dunkl's generalized spherical harmonics}\label{sec3.3}
An alternative approach towards the solution of the eigenvalue problem
for the second order differential
operators $D_1^{\text{C}}$ \eqref{diffA1}, \eqref{diffB1}
is to separate the eigenvalue equation in a `radial' and a `spherical'
component. In essence this is the method used by Calogero \cite{cal:solution}
to obtain the eigenfunctions for the type $A$ Hamiltonian $H_1^{\text{A}}$
\eqref{HA1}
(cf. the comments in the previous subsection).
Specifically, if one substitutes an Ansatz function of the form
$R^{\text{C}}(r) Y^{\text{C}}_l(x)$,
where $R^{\text{C}}(r)$ is a function of $r=\sqrt{x_1^2+\cdots +x_n^2}$ and
$Y^{\text{C}}_l(x)$ is a permutation symmetric homogeneous polynomial of
degree $l$ in $x_1,\ldots ,x_n$ ($C=A$) or
$x_1^2,\ldots ,x_n^2$ ($C=B$), then it is seen that this yields an
eigenfunction
of $D_1^{\text{C}}$ \eqref{diffA1}, \eqref{diffB1} with eigenvalue
$E^{\text{C}}$ if
$R^{\text{C}}(r)$ and $Y^{\text{C}}_l(x)$ satisfy

\begin{subequations}
\begin{eqnarray}\label{radA}
&&-\frac{d^2 R^{\text{A}}}{dr^2} +
\Bigl( 2\omega r -\frac{2l +g_0 n(n-1)+n-1}{r}\Bigr)
\frac{dR^{\text{A}}}{dr} = (E^{\text{A}}-2\omega l) R^{\text{A}} ,\\
\label{radB}
&&-\frac{d^2 R^{\text{B}}}{dr^2} +
\Bigl( 2\omega r -\frac{4l +g_0 n(n-1)+2ng_1+n-1}{r}\Bigr)
\frac{dR^{\text{B}}}{dr} = (E^{\text{B}}-4\omega l) R^{\text{B}}
\end{eqnarray}
\end{subequations}
and
\begin{equation}\label{sphAB}
L^{\text{C}} Y_l^{\text{C}} = 0
\;\;\;\;\;\;\;\;\;\;\;\;\;\;\;\;\;\; (C=A,B)
\end{equation}
with
\begin{subequations}
\begin{eqnarray}\label{sphA}
\makebox[2em]{} L^{\text{A}} \!\!\!\! &=& \!\!\!\! \sum_{1\leq j\leq n}
\frac{\partial^2}{\partial x_j^2}
+2g_0 \sum_{1\leq j< k\leq n} \frac{1}{x_j-x_k}
\Bigl(\frac{\partial}{\partial x_j}-\frac{\partial}{\partial x_k} \Bigr) ,\\
\label{sphB}
\makebox[2em]{} L^{\text{B}}\!\!\!\! &=& \!\!\!\! \sum_{1\leq j\leq n} \left(
\frac{\partial^2}{\partial x_j^2} +2g_1 \frac{\partial}{\partial x_j}\right) \\
&& \!\!\!\!\!\!\!\! +2g_0 \sum_{1\leq j< k\leq n} \left(\frac{1}{x_j-x_k}
\Bigl( \frac{\partial}{\partial x_j}-\frac{\partial}{\partial x_k}\Bigr) +
\frac{1}{x_j+x_k}
\Bigl( \frac{\partial}{\partial x_j}+\frac{\partial}{\partial x_k} \Bigr)
\right) .\nonumber
\end{eqnarray}
\end{subequations}
The `radial' equations \eqref{radA} and \eqref{radB} are confluent
hypergeometric type equations that admit polynomial solutions
for $E^{\text{A}}=2\omega (l +2m)$ and  $E^{\text{B}}=4\omega (l +m)$
(with $m\in \mathbb{N}$) given by Laguerre polynomials in $r^2$
(cf. Eqs.~\eqref{exprepB} and \eqref{lagdiffB})
\begin{subequations}
\begin{eqnarray}
R^{\text{A}}_m(r) &=&
\frac{[l+(n-1)(1+ng_0)/2]_{m}}{(-\omega )^m} \\
&&\times\;
{}_1F_1 \left( \begin{array}{c} -m \\ l + (n-1)(1+ng_0)/2 \end{array}\,
;\, \omega r^2 \right) ,\nonumber \\
R^{\text{B}}_m(r) &=&
\frac{[2l+n(1/2+g_1+(n-1)g_0)]_{m}}{(-\omega )^m} \\
&&\times \; {}_1F_1
\left( \begin{array}{c} -m \\ 2l + n(1/2 +g_1+(n-1)g_0) \end{array}
\, ;\, \omega r^2 \right) .\nonumber
\end{eqnarray}
\end{subequations}
(Here we have chosen the normalization chosen such that $R^{\text{C}}_m(r)$ is
monic.)

The `spherical' equation \eqref{sphAB} was studied
(for type $A$) by Calogero \cite{cal:solution}
and in further detail and more generality (therewith
also including the type $B$) by Dunkl \cite{dun:orthogonal}.
Let $\mathcal{P}_l^{\text{C}}$ be the space of homogeneous
symmetric polynomials of degree $l$ in $x_1,\ldots ,x_n$ ($C=A$) or
$x_1^2,\ldots ,x_n^2$ ($C=B$) and let
$\mathcal{H}_l^{\text{C}}\subset \mathcal{P}_l^{\text{C}}$
be the subspace of polynomials satisfying \eqref{sphAB}.
The polynomials in $\mathcal{H}_l^{\text{C}}$ are referred to as
(symmetric) generalized spherical harmonics. For $g_0,g_1=0$ these generalized
spherical harmonics reduce to ordinary harmonic polynomials in $\mathbb{R}^n$.
It follows from Dunkl's theory in \cite{dun:orthogonal}
(see also \cite{dun:reflection} for the extension to the nonsymmetric case)
that
$\text{dim} (\mathcal{H}_l^{\text{A}})=
\text{dim} (\mathcal{P}_l^{\text{A}})-
\text{dim} (\mathcal{P}_{l-2}^{\text{A}})$
and that
$\text{dim} (\mathcal{H}_l^{\text{B}})=
\text{dim} (\mathcal{P}_l^{\text{B}})-
\text{dim} (\mathcal{P}_{l-1}^{\text{B}})$.
The upshot is that each polynomial $p_\lambda^{\text{C}}$ may be
written uniquely in the form
\begin{subequations}
\begin{eqnarray}\label{decA}
p_\lambda^{\text{A}}(x) &=& \sum_{m=0}^{[ |\lambda |/2 ]}
R_m^{\text{A}}(r) Y^{\text{A}}_{|\lambda|-2m}(x),
\;\;\;\;\;\;\;\;\;\;\;\;\;\;\;\;\;\;\;\;
Y^{\text{A}}_{|\lambda|-2m}\in \mathcal{H}_{|\lambda|-2m}^{\text{A}} \\ [1ex]
\label{decB}
p_\lambda^{\text{B}}(x) &=& \sum_{m=0}^{ |\lambda | }
R_m^{\text{B}}(r) Y^{\text{B}}_{|\lambda|-m}(x),
\;\;\;\;\;\;\;\;\;\;\;\;\;\;\;\;\;\;\;\;\;\;\;
Y^{\text{B}}_{|\lambda|-m}\in \mathcal{H}_{|\lambda|-m}^{\text{B}}
\end{eqnarray}
\end{subequations}
(where $[\cdot ]$ represents the function that extracts the integer part).
Indeed, the functions of the form in the r.h.s. of \eqref{decA}
and \eqref{decB}
are eigenfunctions of $D_1^{\text{C}}$ \eqref{diffA1},
\eqref{diffB1} corresponding to
the eigenvalue $2\omega |\lambda |$ (type $A$) and $4\omega |\lambda |$
(type $B$), respectively. Furthermore, the functions
of this form span a space of dimension
$\sum_{m=0}^{[ |\lambda |/2 ]}
\text{dim}(\mathcal{H}_{|\lambda|-2m}^{\text{A}} )=
\text{dim}(\mathcal{P}^{\text{A}}_{|\lambda |})$
and
$\sum_{m=0}^{ |\lambda | }
\text{dim}(\mathcal{H}_{|\lambda|-m}^{\text{B}} )=
\text{dim}(\mathcal{P}^{\text{B}}_{|\lambda |})$, which is
precisely the multiplicity of the eigenvalues
$E_1^{\text{C}}(\lambda )$ \eqref{evAB1}.

The formulas \eqref{decA}, \eqref{decB} describe the relation between
the Calogero type eigenfunctions of the form $R^{\text{C}}(r)
Y^{\text{C}}_l(x)$
and the Hermite/Laguerre basis $p_\lambda^{\text{C}}(x)$.
To determine precisely which functions $Y^{\text{C}}_l(x)$ appear in the
decompositions
\eqref{decA} and \eqref{decB}, we pick the leading homogeneous parts
on both sides of the equation. It is known from the work of Lassalle
\cite{las:laguerre,las:hermite} that the highest-order homogeneous
part of the multivariable Hermite (type $A$) and Laguerre (type $B$)
polynomials are (with our normalization monic) Jack polynomials
with parameter $\alpha =1/g_0$ in
$x_1,\ldots ,x_n$ and $x_1^2,\ldots ,x_n^2$, respectively.
So, we get upon taking the leading homogeneous part
\begin{subequations}
\begin{eqnarray}\label{declA}
\makebox[2em]{}J_\lambda (x;1/g_0) &=& \sum_{m=0}^{[ |\lambda |/2 ]}
r^{2m} Y^{\text{A}}_{|\lambda|-2m}(x),
\;\;\;\;\;\;\;\;\;\;\;\;\;\;\;\;\;\;\;\;
Y^{\text{A}}_{|\lambda|-2m}\in \mathcal{H}_{|\lambda|-2m}^{\text{A}} \\ [1ex]
\label{declB}
\makebox[2em]{}J_\lambda (x^2;1/g_0) &=& \sum_{m=0}^{ |\lambda | }
r^{2m} Y^{\text{B}}_{|\lambda|-m}(x),
\;\;\;\;\;\;\;\;\;\;\;\;\;\;\;\;\;\;\;\;\;\;\;
Y^{\text{B}}_{|\lambda|-m}\in \mathcal{H}_{|\lambda|-m}^{\text{B}} ,
\end{eqnarray}
\end{subequations}
where $J_\lambda (x;1/g_0)$ and $J_\lambda (x^2;1/g_0)$ denote
the monic Jack polynomial in $x_1,\ldots ,x_n$ and
$x_1^2,\ldots ,x_n^2$ with parameter $\alpha =1/g_0$
\cite{sta:some,mac:symmetric}.
In \cite{dun:orthogonal}, Dunkl
provides inversion formulas for decompositions of the form
\eqref{declA}, \eqref{declB} with which one can express the functions
$Y^{\text{A}}_{|\lambda |-2m}$
and $Y^{\text{B}}_{|\lambda |-m}$ in terms of the homogeneous symmetric
polynomial in the
l.h.s. In our case these inversion formulas become
\begin{subequations}
\begin{eqnarray}\label{DinvA}
Y_{|\lambda |-2m}^{\text{A}}(x) & =&
\Pi_{|\lambda |-2m}^{\text{A}} J_\lambda (x;1/g_0) ,\\
\label{DinvB}
Y_{|\lambda |-m}^{\text{B}}(x) & =&
\Pi_{|\lambda |-m}^{\text{B}} J_\lambda (x^2;1/g_0)
\end{eqnarray}
\end{subequations}
with
\begin{eqnarray*}
\Pi_{|\lambda |-2m}^{\text{A}} &=&
\frac{1}{4^m m!\: [n/2 +d^{\text{A}} +|\lambda |-2m]_m}
T_{|\lambda |-2m}^{\text{A}} (L^{\text{A}})^m \\
T_{k}^{\text{A}} &=&
\sum_{j=0}^{[k/2]} \frac{r^{2j}}{4^j j! \: [-n/2 -d^{\text{A}} -k+2]_j}
(L^{\text{A}})^j ,
\end{eqnarray*}
$d^{\text{A}}= g_0 n(n-1)/2 $ and
\begin{eqnarray*}
\Pi_{|\lambda |-m}^{\text{B}} &=&
\frac{1}{4^m m!\: [n/2 +d^{\text{B}} +2|\lambda |-2m]_m}
T_{|\lambda |-m}^{\text{B}} (L^{\text{B}})^m \\
T_{k}^{\text{A}} &=&
\sum_{j=0}^{k} \frac{r^{2j}}{4^j j!\: [-n/2 -d^{\text{B}} -2k+2]_j}
(L^{\text{B}})^j ,
\end{eqnarray*}
$d^{\text{B}}= g_0 n(n-1) +ng_1 $.
Formulas \eqref{decA}, \eqref{decB} combined with \eqref{DinvA}, \eqref{DinvB}
render the decomposition of the multivariable Hermite and Laguerre polynomials
in terms of Dunkl's generalized spherical harmonics in a closed form.

\subsection{Notes}\label{sec3.4}
{\em i. Orthogonality.}
The orthogonality of the multivariable Hermite and Laguerre polynomials
was stated in \cite{mac:hypergeometric} and \cite{las:laguerre,las:hermite}.
It was proved by Macdonald for the cases $g_0=1/2$, $1$ and $2$.
Recently, Baker and Forrester \cite{bak-for:calogero-sutherland} proposed a
proof
valid for general parameters that exploits
the fact that the polynomials may be seen as limiting cases of the
multivariable Jacobi polynomials
\cite{vre:formulas,deb:systeme,bee-opd:certain,hec:elementary}.
In essence this approach should boil down to extending the transition in
\cite{mac:some}
from the Selberg integral (whose integrand is the weight function for
the Jacobi polynomials) to the Mehta-Macdonald integrals \eqref{mehtaA},
\eqref{mehtaB}
(whose integrands are the weight functions for the multivariable
Hermite and Laguerre polynomials) so as to include also the polynomials of
higher degree.

{\em ii. Norm formulas.}
Norm formulas for the multivariable Hermite and Laguerre polynomials
can be found in \cite{mac:hypergeometric} and \cite{las:laguerre,las:hermite},
with again proofs given by Macdonald for $g_0=1/2$, $1$ and $2$.
Baker and Forrester provide a proof of these norm formulas
for general parameters that hinges on a generating function approach
(and uses also the orthogonality) \cite{bak-for:calogero-sutherland}.
The expressions for the norms given in
\cite{mac:hypergeometric,las:laguerre,las:hermite,bak-for:calogero-sutherland}
are written in terms of the Mehta-Macdonald integrals \eqref{mehtaA},
\eqref{mehtaB}
and evaluations of Jack symmetric functions at the identity.
It can be verified with the aid of
known evaluation formulas for the Jack symmetric functions
due to Stanley (see \cite{sta:some,mac:symmetric}) that
our norm formulas in Theorem~\ref{normthm} are in agreement with those
in
\cite{mac:hypergeometric,las:laguerre,las:hermite,bak-for:calogero-sutherland}.

{\em iii. Differential equations.}
The second order differential equation for the multivariable Hermite and
Laguerre
polynomials of the form
$D_1^{\text{C}} p_\lambda^{\text{C}}= E_1^{\text{C}}(\lambda )
p_\lambda^{\text{C}}$
was already given by Macdonald \cite{mac:hypergeometric} and Lassalle
\cite{las:laguerre,las:hermite}.
In \cite{bak-for:calogero-sutherland} a procedure is described to obtain the
complete system
of $n$-independent differential equations starting from the explicitly known
differential equations for the Jack symmetric functions found by
Sekiguchi \cite{sek:zonal} and Debiard \cite{deb:polynomes}
(see also \cite{mac:commuting,mac:symmetric}).
Another approach to arrive at such a system of differential equations for
$p_\lambda^{\text{C}}$
is to employ Dunkl's differential-reflection operators; see \cite{kak:common}
for a treatment along this lines of the Hermite case.
Yet another strategy would be
to analyze the limit behavior of the system of differential equations for the
Jacobi
polynomials due to Debiard \cite{deb:systeme}, with respect to the transitions
`Jacobi' $\rightarrow$ `Hermite' and `Jacobi' $\rightarrow$ `Laguerre'.
All these methods have in common
with our approach in Section~\ref{sec2.2} that it seems a priori
difficult (from a computational point of view)
to extract further explicit information pertaining to a possible closed form
for the
higher-order differential equations at the confluent hypergeometric level
(cf. also the comments in the last paragraph of Section~\ref{sec3.2}).

{\em iv. Recurrence formulas.}
Recurrence relations of the type given by Theorem~\ref{recrthm1} (i.e. the
simplest
ones, corresponding to the first elementary symmetric function) were
recently derived independently by Baker and Forrester using
a generating function for the polynomials \cite{bak-for:calogero-sutherland}.
As it stands, their recurrence formulas are a little less explicit
than those obtained here since the coefficients
are written in terms of certain implicitly defined generalized binomial
coefficients
and furthermore contain evaluations of the Jack symmetric
functions at the identity.
In order
make their formulas fully explicit (so as to compare with
Theorem~\ref{recrthm1}) one again needs Stanley's expressions
for the Jack symmetric functions at the identity in combination with an
explicit representation for the specific binomial coefficients at hand that can
be found in
\cite{las:formule}.

{\em v. Rodrigues formulas.} Recently, Ujino and Wadati derived Rodrigues type
formulas for
the multivariable Hermite polynomials \cite{uji-wad:rodrigues}
following an approach due to Lapointe and Vinet who obtained similar Rodrigues
formulas for the Jack symmetric functions
\cite{lap-vin:rodrigues,lap-vin:exact}.
Such Rodrigues formulas are particularly useful when trying to answer questions
regarding the structure
of the coefficients $c_{\lambda ,\mu}$  that appear in the expansion of
the polynomials in terms of monomial symmetric functions.
For instance, the Rodrigues formulas allowed Lapointe and Vinet
to prove a weak form of
the Macdonald-Stanley conjecture saying that (in an appropriate normalization)
the
expansion coefficients for the Jack symmetric functions in terms of monomial
symmetric functions are
polynomials in the parameters with integer coefficients
\cite{lap-vin:rodrigues}.
(See \cite{sta:some,mac:symmetric} for the Macdonald-Stanley conjecture and
various related conjectures.)
A similar statement also holds true for the multivariable Hermite polynomials
\cite{uji-wad:rodrigues}.
\vfill

\section{Proofs}\label{sec4}
In this section the properties of the multivariable confluent hypergeometric
families stated in Section~\ref{sec2} are proven by viewing the polynomials
as degenerate (limiting) cases of the multivariable hypergeometric
continuous Hahn families (type $A$) and Wilson families (type $B$)
that were investigated in \cite{die:multivariable,die:properties}.

\subsection{Orthogonality properties}\label{sec4.1}
In \cite{die:multivariable,die:properties} multivariable continuous Hahn
and Wilson type polynomials were considered that are
associated to the weight functions
\begin{subequations}
\begin{equation}\label{weightcH}
\Delta^{\text{cH}}(x) =
\prod_{1\leq j< k\leq n}
\left| \frac{\Gamma (g_0 + i(x_j-x_k))}{\Gamma (i(x_j-x_k))} \right|^2 \;
\prod_{1\leq j\leq n}
\left| \Gamma (a+ix_j)\: \Gamma (b+ix_j ) \right|^2
\end{equation}
and
\begin{eqnarray}\label{weightW}
\Delta^{\text{W}}(x)\!\!\! &=&\!\!\!
\prod_{1\leq j< k\leq n}
\left| \frac{\Gamma (g_0 + i(x_j-x_k))\: \Gamma (g_0 + i(x_j+x_k))}
            {\Gamma (i(x_j-x_k))\: \Gamma (i(x_j+x_k))} \right|^2 \\
& & \makebox[1em]{} \times \prod_{1\leq j\leq n}
\left| \frac{\Gamma (a+ix_j)\, \Gamma (b+ix_j )\,
             \Gamma (c+ix_j)\, \Gamma (d+ix_j )}{\Gamma (2ix_j)} \right|^2
\nonumber
\end{eqnarray}
\end{subequations}
(with $g_0\geq 0$ and $\text{Re}(a,b,c,d)> 0$).
Specifically, the multivariable continuous Hahn polynomials are defined by the
conditions
$A.1$, $A.2$ in Section~\ref{sec2.1} with
$\Delta^{\text{cH}}(x)$ \eqref{weightcH} replacing the weight function
$\Delta^{\text{A}}(x)$ \eqref{weightA}. Similarly, the multivariable
Wilson polynomials are defined by the conditions $B.1$, $B.2$ in
Section~\ref{sec2.1} with $\Delta^{\text{B}}(x)$ \eqref{weightB} being
replaced by $\Delta^{\text{W}}(x)$ \eqref{weightW}.

If we rescale the variables by substituting
\begin{equation}\label{scale}
x_j \longrightarrow x_j/\beta ,\;\;\;\;\;\;\;\;\;\;\;\;\;\; j=1,\ldots ,n
\end{equation}
and simultaneously perform a reparametrization of the form
\begin{equation}\label{repa}
a = (\beta^{2} \varpi )^{-1},\;\;\;\;\;
b = (\beta^{2} \varpi^\prime )^{-1},\;\;\;\;\;
c= \text{g}_1,\;\;\;\;\;
d= \text{g}_1^\prime +1/2
\end{equation}
(with $\varpi, \varpi^\prime >0$, $\text{g}_1, \text{g}_1^\prime \geq 0$ and
$\beta$ real),
then the weight functions $\Delta^{\text{cH}}(x)$ \eqref{weightcH} and
$\Delta^{\text{W}}(x)$ \eqref{weightW}
pass (upon multiplication by the overall
normalization constants $\mathcal{D}^{\text{A}}(\beta )$ and
$\mathcal{D}^{\text{B}}(\beta )$) over
into
\begin{subequations}
\begin{eqnarray}\label{weightAbeta}
\Delta^{\text{A}}_\beta(x)&=& \mathcal{D}^{\text{A}} (\beta)
\prod_{1\leq j< k\leq n}
\left| \frac{\Gamma (g_0 + i\beta^{-1}(x_j-x_k) )}{\Gamma (i\beta^{-1}(x_j-x_k)
)} \right|^2 \\
 && \makebox[1em]{} \times \prod_{1\leq j\leq n}
           \left| \Gamma (\frac{1}{\varpi \beta^2}+i\beta^{-1}x_j)\:
                  \Gamma (\frac{1}{\varpi^\prime \beta^2}+i\beta^{-1} x_j )
\right|^2 \nonumber
\end{eqnarray}
and
\begin{eqnarray}
\Delta^{\text{B}}_\beta (x)&=& \mathcal{D}^{\text{B}}(\beta)
\prod_{1\leq j< k\leq n}
\left| \frac{\Gamma (g_0 + i\beta^{-1}(x_j-x_k))\: \Gamma (g_0 +
i\beta^{-1}(x_j+x_k))}
            {\Gamma (i\beta^{-1}(x_j-x_k))\: \Gamma (i\beta^{-1}(x_j+x_k))}
\right|^2 \nonumber \\
&& \makebox[1em]{} \times \prod_{1\leq j\leq n}
\left| \frac{\Gamma (\text{g}_1 +i\beta^{-1}x_j)\, \Gamma (\text{g}_1^\prime
+1/2 +i\beta^{-1}x_j )}
            {\Gamma (i\beta^{-1} x_j)\, \Gamma (1/2 +i\beta^{-1}x_j) } \right.
\label{weightBbeta} \\
&& \makebox[7em]{}  \times
\left.   \Gamma (\frac{1}{\varpi \beta^2}+i\beta^{-1}x_j)\:
         \Gamma (\frac{1}{\varpi^\prime \beta^2}+i\beta^{-1} x_j ) \right|^2  ,
\nonumber
\end{eqnarray}
\end{subequations}
respectively.
The normalization constants $\mathcal{D}^{\text{A}}(\beta ) $ and
$\mathcal{D}^{\text{B}} (\beta )$
are introduced so as to ensure finite limiting behavior
of the
weight functions for $\beta \rightarrow 0$.
It is not so difficult to check (see appendix)---using
Stirling's formula for the asymptotics
of $\Gamma (z)$ for $|z|\rightarrow \infty$
(see e.g. \cite{abr-ste:handbook})---that if one takes
\begin{eqnarray*}
\mathcal{D}^{\text{A}}(\beta )&=&
|\beta |^{g_0 n(n-1)}
\delta (\varpi ,\beta  )^{2n} \delta (\varpi^\prime ,\beta )^{2n} ,\\
\mathcal{D}^{\text{B}}(\beta )&=&
|\beta |^{2 g_0n(n-1) +2n(\text{g}_1+\text{g}_1^\prime )}
          \delta (\varpi ,\beta )^{2n} \delta (\varpi^\prime ,\beta  )^{2n}
\end{eqnarray*}
where
\begin{equation}
\delta (\alpha , \beta )
\equiv
\sqrt{\frac{e}{2\pi}}
e^{(1+\log (\beta^2\alpha ) ) ( 1/(\beta^{2}\alpha ) - 1/2)},
\end{equation}
then the weight functions
$\Delta^{\text{A}}_\beta (x)$ \eqref{weightAbeta}
and $\Delta^{\text{B}}_\beta (x)$ \eqref{weightBbeta} converge pointwise to
$\Delta^{\text{A}}(x)$ \eqref{weightA} and $\Delta^{\text{B}}$ \eqref{weightB}
for $\beta\rightarrow 0$
provided the parameters of
$\Delta^{\text{C}}(x)$ are related
to those of $\Delta^{\text{C}}_\beta (x)$ by
\begin{equation}\label{corpa}
\omega\equiv \varpi + \varpi^\prime \;\;\;\;\;\;\;\;\;\; \text{and}
\;\;\;\;\;\;\;\;\;\;
g_1 \equiv \text{g}_1 + \text{g}_1^\prime   .
\end{equation}
For our purposes, however, pointwise convergence is not sufficient
and we need a somewhat stronger convergence result stating that the
corresponding
measures pass over into each other:
\begin{eqnarray}\label{convmeas}
& &\int_{-\infty}^{\infty} \cdots \int_{-\infty}^{\infty}
p(x)\: \Delta^{\text{C}}(x)\: dx_1\cdots dx_n
 =  \\
&& \makebox[2em]{} \lim_{\beta\rightarrow 0}\;
\int_{-\infty}^{\infty} \cdots \int_{-\infty}^{\infty}
p(x)\: \Delta^{\text{C}}_\beta(x)\: dx_1\cdots dx_n
\;\;\;\;\;\; (C= A\; \text{or}\; B),  \nonumber
\end{eqnarray}
where $p(x)$ denotes an arbitrary polynomial in the variables $x_1,\ldots ,
x_n$.
A proof of the limit formula \eqref{convmeas} can be found in the appendix at
the end of
the paper.

Now, let $\{ p_{\lambda ,\beta}^{\text{A}} \}_{\lambda\in\Lambda}$ and
$\{ p_{\lambda ,\beta}^{\text{B}} \}_{\lambda\in\Lambda}$
be the bases determined by the conditions A.1, A.2 and B.1, B.2 in
Section~\ref{sec2.1}
with the weight functions
$\Delta^{\text{A}}(x)$ \eqref{weightA} and $\Delta^{\text{B}}$ \eqref{weightB}
being replaced by
$\Delta^{\text{A}}_\beta (x)$ \eqref{weightAbeta}
and $\Delta^{\text{B}}_\beta (x)$ \eqref{weightBbeta}, respectively.
So, up to scaling and reparametrization,
the polynomials $p_{\lambda ,\beta}^{\text{A}}(x)$
amount to the multivariable continuous Hahn polynomials (multiplied
by $\beta^{|\lambda |}$) and
the polynomials $p_{\lambda ,\beta}^{\text{B}}(x)$ amount to the multivariable
Wilson
polynomials (multiplied by $\beta^{2|\lambda |}$) from
\cite{die:multivariable,die:properties}.
It then follows from the defining properties for the polynomials
(of the form A.1, A.2 and B.1, B.2) and the limit formula \eqref{convmeas} that
\begin{equation}\label{pollim}
p_{\lambda }^{\text{C}} (x) = \lim_{\beta \rightarrow 0}\; p_{\lambda
,\beta}^{\text{C}} (x)
\end{equation}
and that
\begin{eqnarray} \label{iplim}
\langle p_\lambda^{\text{C}} ,  p_{\mu}^{\text{C}} \rangle_{\text{C}}  &=&
\lim_{\beta \rightarrow 0}\;
\langle p_{\lambda ,\beta}^{\text{C}} ,  p_{\mu ,\beta}^{\text{C}}
\rangle_{\text{C},\beta} \\
&=&
\lim_{\beta \rightarrow 0}\;
\int_{-\infty}^{\infty} \cdots \int_{-\infty}^{\infty}
p_{\lambda ,\beta}^{\text{C}}(x)\:
\overline{p_{\mu ,\beta }^{\text{C}} (x)}\: \Delta^{\text{C}}_\beta (x)\:
dx_1\cdots dx_n \nonumber
\end{eqnarray}
(of course again assuming a correspondence between the parameters in accordance
with \eqref{corpa}).
Theorem~\ref{orthothm} is now immediate from limit formula \eqref{iplim}
and the orthogonality of
the multivariable continuous Hahn and Wilson families \cite{die:properties}
(which translates into the orthogonality of the bases
$\{ p_{\lambda ,\beta}^{\text{A}} \}_{\lambda\in\Lambda}$ and
$\{ p_{\lambda ,\beta}^{\text{B}} \}_{\lambda\in\Lambda}$ with respect to
the weight functions $\Delta^{\text{A}}_\beta (x)$ \eqref{weightAbeta}
and $\Delta^{\text{B}}_\beta (x)$ \eqref{weightBbeta}).

The proof of Theorem~\ref{normthm} is based on explicit expressions
for the quotients
$\langle p_{\lambda }^{\text{cH}} ,  p_{\lambda }^{\text{cH}}
\rangle_{\text{cH}}/
\langle 1 ,  1 \rangle_{\text{cH}}$ and
$\langle p_{\lambda }^{\text{W}} ,  p_{\lambda }^{\text{W}} \rangle_{\text{W}}/
\langle 1 ,  1 \rangle_{\text{W}}$---i.e. the
ratios of the squared norm of the multivariable continuous Hahn resp. Wilson
polynomials
and the unit polynomial (with respect to the $L^2$ inner product over
$\mathbb{R}^n$
with weight function $\Delta^{\text{cH}}(x)$ resp.
$\Delta^{\text{W}}(x)$)---that
were computed in \cite{die:properties}.
We have
\begin{eqnarray*}
\frac{\langle p_{\lambda }^{\text{cH}} ,  p_{\lambda }^{\text{cH}}
\rangle_{\text{cH}}}
     {\langle 1 ,  1 \rangle_{\text{cH}}}
\!\!\!\! &=&\!\!\!\! 2^{-4|\lambda |}
\prod_{1\leq j < k\leq n} \left(
\frac{[ g_0 +\rho_j^{\text{cH}}+\rho_k^{\text{cH}},
         1-g_0+\rho_j^{\text{cH}}+\rho_k^{\text{cH}}]_{\lambda_j+\lambda_k}}
     {[ \rho_j^{\text{cH}}+\rho_k^{\text{cH}},
     1+\rho_j^{\text{cH}}+\rho_k^{\text{cH}}]_{\lambda_j+\lambda_k}} \right. \\
&& \makebox[5em]{} \times \left.
\frac{[ g_0 +\rho_j^{\text{cH}}-\rho_k^{\text{cH}},
      1-g_0+\rho_j^{\text{cH}}-\rho_k^{\text{cH}}]_{\lambda_j-\lambda_k}}
     {[ \rho_j^{\text{cH}}-\rho_k^{\text{cH}},
        1+\rho_j^{\text{cH}}-\rho_k^{\text{cH}}]_{\lambda_j-\lambda_k}} \right)
\\
&&  \!\!\!\!\!\!\!\! \times \prod_{1\leq j\leq n}
\frac{ [\hat{a}^{\text{cH}}+\rho_j^{\text{cH}},
        \hat{b}^{\text{cH}}+\rho_j^{\text{cH}},
       1-\hat{a}^{\text{cH}}+\rho_j^{\text{cH}} ,
       1-\hat{b}^{\text{cH}}+\rho_j^{\text{cH}}]_{\lambda_j} }
     { [\rho_j^{\text{cH}} ,1+  \rho_j^{\text{cH}} ]_{\lambda_j}}
\end{eqnarray*}
with $\rho_j^{\text{cH}}=(n-j)g_0+a+b-1/2$,
$\hat{a}^{\text{cH}}=a+b-1/2$ and $\hat{b}^{\text{cH}}=a-b+1/2$; and we
have
\begin{eqnarray*}
\frac{\langle p_{\lambda }^{\text{W}} ,  p_{\lambda }^{\text{W}}
\rangle_{\text{W}}}
     {\langle 1 ,  1 \rangle_{\text{W}}} \!\!\!\! &=& \!\!\!\!
\prod_{1\leq j < k\leq n} \left(
\frac{[ g_0 +\rho_j^{\text{W}}+\rho_k^{\text{W}},
        1-g_0+\rho_j^{\text{W}}+\rho_k^{\text{W}}]_{\lambda_j+\lambda_k}}
     {[ \rho_j^{\text{W}}+\rho_k^{\text{W}},
        1+\rho_j^{\text{W}}+\rho_k^{\text{W}}]_{\lambda_j+\lambda_k}} \right.
\\
&& \makebox[5em]{} \times \left.
\frac{[ g_0 +\rho_j^{\text{W}}-\rho_k^{\text{W}},
       1-g_0+\rho_j^{\text{W}}-\rho_k^{\text{W}}]_{\lambda_j-\lambda_k}}
     {[ \rho_j^{\text{W}}-\rho_k^{\text{W}},
      1+\rho_j^{\text{W}}-\rho_k^{\text{W}}]_{\lambda_j-\lambda_k}} \right) \\
&& \!\!\!\!\!\!\!\! \times \prod_{1\leq j\leq n}
\left( \frac{ [\hat{a}^{\text{W}}+\rho_j^{\text{W}},
               \hat{b}^{\text{W}}+\rho_j^{\text{W}},
               \hat{c}^{\text{W}}+\rho_j^{\text{W}},
               \hat{d}^{\text{W}}+\rho_j^{\text{W}}]_{\lambda_j} }
     { [2\rho_j^{\text{W}} ]_{2\lambda_j}} \right. \\
&& \makebox[1em]{} \times \left.
\frac{ [1-\hat{a}^{\text{W}}+\rho_j^{\text{W}},
        1-\hat{b}^{\text{W}}+\rho_j^{\text{W}},
        1-\hat{c}^{\text{W}}+\rho_j^{\text{W}},
        1-\hat{d}^{\text{W}}+\rho_j^{\text{W}}]_{\lambda_j} }
     { [1+2\rho_j^{\text{W}} ]_{2\lambda_j}} \right)
\end{eqnarray*}
with $\rho_j^{\text{W}}=(n-j)g_0 + (a+b+c+d-1)/2$,
$\hat{a}^{\text{W}}= (a-b+c-d+1)/2$,
$\hat{b}^{\text{W}}=(-a+b+c-d+1)/2$,
$\hat{c}^{\text{W}}=(a+b+c+d-1)/2$
and $\hat{d}^{\text{W}}=(-a-b+c+d+1)/2$.
Here we have adopted the notation
\begin{equation}
 [a_1,\ldots ,a_p]_l\equiv [a_1]_l \cdots [a_p]_l
\end{equation}
($[a]_0= 0$, $[a]_l= a (a+1)\cdots (a+l-1)$)
for the Pochhammer symbols (and as usual
$|\lambda | \equiv \lambda_1+\cdots +\lambda_n$).
Scaling of the variables and reparametrization in accordance with \eqref{scale}
and \eqref{repa}
leads us to the corresponding expressions for $p_{\lambda ,\beta }^{\text{A}}$
and
$p_{\lambda ,\beta }^{\text{B}}$, which
entail for $\beta \rightarrow 0$ (using \eqref{iplim} and \eqref{corpa})
\begin{subequations}
\begin{eqnarray}\label{ratA}
\makebox[1em]{}
\frac{\langle p_{\lambda }^{\text{A}} ,  p_{\lambda }^{\text{A}}
\rangle_{\text{A}}}
     {\langle 1 ,  1 \rangle_{\text{A} }}\!\!\!\! &=&\!\!\!\!
(2\omega)^{-|\lambda |}
\prod_{1\leq j <k\leq n}
\frac{ [ (k-j+1)g_0, 1+(k-j-1)g_0 ]_{\lambda_j -\lambda_k} }
     { [ (k-j)g_0, 1+ (k-j)g_0 ]_{\lambda_j -\lambda_k}    } \\
& & \times \prod_{1\leq j\leq n} [ 1+ (n-j)g_0 ]_{\lambda_j} \nonumber
\end{eqnarray}
and
\begin{eqnarray}\label{ratB}
\makebox[1em]{}
\frac{\langle p_{\lambda }^{\text{B}} ,  p_{\lambda }^{\text{B}}
\rangle_{\text{B}}}
     {\langle 1 ,  1 \rangle_{\text{B} }}\!\!\!\! &=&\!\!\!\!
(2\omega)^{-2|\lambda |}
\prod_{1\leq j <k\leq n}
\frac{ [ (k-j+1)g_0, 1+(k-j-1)g_0 ]_{\lambda_j -\lambda_k} }
     { [ (k-j)g_0, 1+ (k-j)g_0 ]_{\lambda_j -\lambda_k}    } \\
& & \times \prod_{1\leq j\leq n} [ (n-j)g_0 +g_1+1/2, 1+ (n-j)g_0 ]_{\lambda_j}
. \nonumber
\end{eqnarray}
\end{subequations}
Combination of
the expressions for the ratios in \eqref{ratA} and \eqref{ratB} with
the Mehta-Macdonald
formulas for $\langle 1 ,  1 \rangle_{\text{A} }$ and
$\langle 1 ,  1 \rangle_{\text{B} }$ (cf. \eqref{mehtaA}, \eqref{mehtaB}),
produces
evaluation formulas for $\langle p_{\lambda }^{\text{A}} ,  p_{\lambda
}^{\text{A}} \rangle_{\text{A}}$
and $\langle p_{\lambda }^{\text{B}} ,  p_{\lambda }^{\text{B}}
\rangle_{\text{B}}$ that
can be cast in the form given by Theorem~\ref{normthm}.
Indeed, one easily checks that
the norm formulas in Theorem~\ref{normthm} are in agreement with the ratio
formulas
\eqref{ratA}, \eqref{ratB} (using the relation $\Gamma (a+l)/\Gamma (a)=[a]_l$)
and, furthermore, that they boil down to the Mehta-Macdonald
evaluation formulas \eqref{mehtaA}, \eqref{mehtaB} for $\lambda =\mathbf{0}$.
In order
to verify the reduction to \eqref{mehtaA} and \eqref{mehtaB} for
$\lambda =\mathbf{0}$, one uses the identity
\begin{equation*}
n! \prod_{1\leq j < k\leq n}
    \frac{\Gamma ( (k-j+1)g_0 )\:
          \Gamma ( 1+ (k-j-1)g_0 )}
         {\Gamma ( (k-j)g_0 )\:
          \Gamma ( 1+ (k-j)g_0 )} = \frac{\Gamma (1+ng_0)}{(\Gamma (1+g_0) )^n}
\end{equation*}
(which is derived by
canceling common factors in the numerator and denominator of the l.h.s.
and some further manipulations involving the standard shift property
for the gamma function $z\, \Gamma (z) = \Gamma (z+1)$).

\subsection{Differential equations}\label{sec4.2}
In \cite{die:properties} systems of difference equations for
the multivariable continuous Hahn and Wilson polynomials associated
to the weight functions $\Delta^{\text{cH}}(x)$ \eqref{weightcH}
and $\Delta^{\text{W}}(x)$ \eqref{weightW} were presented. If we
rescale the variables and reparametrize in accordance with
\eqref{scale} and \eqref{repa}, then the difference equations in
question pass
over into difference equations of the form
\begin{equation}\label{diffeqbeta}
D^{\text{C}}_{r ,\beta} p_{\lambda ,\beta}^{\text{C}} =
E^{\text{C}}_{r,\beta}(\lambda) p_{\lambda  ,\beta}^{\text{C}}
\;\;\;\;\;\;\;\;\;\;\;\;\;\; (r=1,\ldots ,n;\; C=A,B),
\end{equation}
for the polynomials $p_{\lambda ,\beta}^{\text{A}}(x)$ and
$p_{\lambda ,\beta}^{\text{B}}(x)$
associated to the weight functions
$\Delta^{\text{A}}_\beta (x)$ \eqref{weightAbeta} and
$\Delta^{\text{B}}_\beta (x)$ \eqref{weightBbeta}.
Here the difference operators $D^{\text{A}}_{r ,\beta}$ and
$D^{\text{B}}_{r ,\beta}$ in the l.h.s. are the same as in Section~\ref{sec2.2}
but
with $w^{\text{A}}(z)$ and
$w^{\text{B}}(z)$ replaced by
\begin{subequations}
\begin{equation}\label{wArep}
w^{\text{A}}(z) =
 (1+i\beta\varpi z) (1+i\beta\varpi^\prime z)
\end{equation}
and
\begin{equation}\label{wBrep}
w^{\text{B}}(z) =
\Bigl( 1 +\frac{\beta \text{g}_1}{iz}\Bigr)
\Bigl( 1 +\frac{\beta \text{g}_1^\prime }{(iz+\beta /2)}\Bigr)
(1+i\beta\varpi z) (1+i\beta\varpi^\prime z)  .
\end{equation}
\end{subequations}
The corresponding eigenvalues in the r.h.s. are given by
\begin{subequations}
\begin{eqnarray}\label{evAbeta}
\makebox[1em]{} E^{\text{A}}_{r,\beta}(\lambda)\!\!\!\! &=&\!\!\!\!
(\beta^4 \varpi \varpi^\prime)^r
E_r((\rho^{\text{A}}_1+\lambda_1)^2, \ldots ,(\rho^{\text{A}}_n+\lambda_n)^2 ;
    (\rho^{\text{A}}_r)^2,\ldots ,(\rho^{\text{A}}_n)^2 ) ,\\
\makebox[1em]{} E^{\text{B}}_{r,\beta}(\lambda)\!\!\!\! &=&\!\!\!\!
(4\beta^4 \varpi \varpi^\prime)^r
E_r((\rho^{\text{B}}_1+\lambda_1)^2, \ldots ,(\rho^{\text{B}}_n+\lambda_n)^2 ;
    (\rho^{\text{B}}_r)^2,\ldots ,(\rho^{\text{B}}_n)^2 ) \label{evBbeta}
\end{eqnarray}
\end{subequations}
with
\begin{eqnarray}\label{Er}
\lefteqn{E_r(\zeta_1,\ldots ,\zeta_n;\eta_r,\ldots ,\eta_n) \equiv } && \\
&& \sum_{\stackrel{J\subset \{ 1,\ldots ,n\}}{0\leq |J|\leq r}}
(-1)^{r-|J|}
\prod_{j\in J}\zeta_j \sum_{r\leq l_1\leq\cdots\leq l_{r-|J|}\leq n}
\eta_{l_1}\cdots \eta_{l_{r-|J|}} \nonumber
\end{eqnarray}
and
\begin{subequations}
\begin{eqnarray}
\label{rhoA}
\rho^{\text{A}}_j&=& (n-j)g_0 + \beta^{-2}\Bigl( \frac{1}{\varpi} +
 \frac{1}{\varpi^\prime} \Bigr) -1/2, \\
\label{rhoB}
\rho^{\text{B}}_j&=& (n-j)g_0 + (\text{g}_0+\text{g}_0^\prime )/2 +
\beta^{-2}\Bigl( \frac{1}{2\varpi} +
 \frac{1}{2\varpi^\prime} \Bigr) -1/4 .
\end{eqnarray}
\end{subequations}

In order to derive the differential equations
for the multivariable Hermite and Laguerre polynomials
$p_\lambda^{\text{A}}$ and $p_\lambda^{\text{B}}$, we
shall analyze the behavior of the difference
equations for $p_{\lambda ,\beta}^{\text{A}}$ and
$p_{\lambda ,\beta}^{\text{B}}$ when the step size parameter
$\beta$ tends to zero.
Essential is the behavior
of the eigenvalues $E^{\text{C}}_{r,\beta}(\lambda )$
\eqref{evAbeta}, \eqref{evBbeta}, which is governed by the
following lemma.
\begin{lemma}\label{limevlem}
One has
\begin{equation*}
\lim_{\beta \rightarrow 0} \beta^{-2r} E^{\text{C}}_{r,\beta}(\lambda ) =
E^{\text{C}}_{r}(\lambda )
\;\;\;\;\;\;\;\;\;\;\;\;\;\;\;\; (r=1,\ldots ,n;\; C=A,B),
\end{equation*}
where $E^{\text{C}}_{r}(\lambda )$ is
given by the expression in Theorem~\ref{diffeqthm} with
$\omega = \varpi +\varpi^\prime$.
\end{lemma}
\begin{proof}
The proof hinges on the following decomposition of $E_r(\cdots ;\cdots )$
\eqref{Er} (cf. \cite[Lemma B.2]{die:commuting})
\begin{eqnarray*}
\lefteqn{E_r(\zeta_1,\ldots ,\zeta_n;\eta_r,\ldots ,\eta_n) = } && \\
&& (\zeta_n-\eta_n) \times\Bigl(
\sum_{\stackrel{J\subset \{ 1,\ldots ,n-1\}}{0\leq |J|\leq r-1}}
(-1)^{r-1-|J|}
\prod_{j\in J}\zeta_j \sum_{r\leq l_1\leq\cdots\leq l_{r-1-|J|}\leq n}
\eta_{l_1}\cdots \eta_{l_{r-1-|J|}} \Bigr) \nonumber \\
&&\makebox[3em]{} + \sum_{\stackrel{J\subset \{ 1,\ldots ,n-1\}}{0\leq |J|\leq
r}}
(-1)^{r-|J|}
\prod_{j\in J}\zeta_j \sum_{r\leq l_1\leq\cdots\leq l_{r-|J|}\leq n-1}
\eta_{l_1}\cdots \eta_{l_{r-|J|}} ,\nonumber
\end{eqnarray*}
with the convention that the sum in the second line equals one if $r=1$ and
that
the sum in the third line equals zero if $r=n$. Notice that these sums
are of the same type as in the r.h.s. of \eqref{Er} but with one
$\zeta$-variable less
(there is no longer dependence on $\zeta_n$).
With the aid of this formula and induction on the number of variables it is
not difficult to infer that
\begin{eqnarray*}
\lim_{\beta\rightarrow 0}
\frac{E^{\text{C}}_{r ,\beta} (\lambda )}{\beta^{2r} (\varpi+ \varpi^\prime)^r
M^r}
&=& \lambda_n\times \Bigl(
\sum_{\stackrel{J\subset \{1 ,\ldots ,n-1\}}{|J|=r-1}}
         \prod_{j\in J} \lambda_j \Bigr)
    +\sum_{\stackrel{J\subset \{1 ,\ldots ,n-1\}}{|J|=r}} \prod_{j\in J}
\lambda_j \\
&=& \sum_{\stackrel{J\subset \{1 ,\ldots ,n\}}{|J|=r}} \prod_{j\in J} \lambda_j
,
\end{eqnarray*}
where $M=2$ for $C=A$ and $M=4$ for $C=B$.
\end{proof}

It is immediate from Lemma~\ref{limevlem} and
the limit formula \eqref{pollim} that
$E^{\text{C}}_{r,\beta}(\lambda ) p_{\lambda ,\beta}^{\text{C}}=
\beta^{2r} E^{\text{C}}_{r}(\lambda ) p_{\lambda}^{\text{C}}
\; + o(\beta^{2r})$ (assuming of course the usual identification
of the parameters via $\omega =\varpi +\varpi^\prime$
and $g_1=\text{g}_1+\text{g}_1^\prime$).
If we plug this asymptotics
into the eigenvalue equation \eqref{diffeqbeta}, then we arrive
at the following asymptotic behavior of
the corresponding l.h.s. for $\beta\rightarrow 0$
\begin{equation}\label{limdiffbasis}
(D^{\text{C}}_{r,\beta } p_{\lambda ,\beta}^{\text{C}})(x) =
\beta^{2r} (D^{\text{C}}_{r } p_{\lambda }^{\text{C}})(x)
\;\;\; + o(\beta^{2r}),
\end{equation}
with $D^{\text{C}}_{r }$ a certain partial differential operator.
In other words, for $\beta\rightarrow 0$
the difference equation for $p_{\lambda ,\beta}^{\text{C}}$
passes over (cf. Corollary~\ref{limdiffcor}) into a differential equation for
the polynomial $p_{\lambda }^{\text{C}}$ of the form
$D^{\text{C}}_{r }p_{\lambda }^{\text{C}}=
E^{\text{C}}_{r}(\lambda )p_{\lambda }^{\text{C}}$, where
$D^{\text{C}}_{r }$ is a certain differential operator and
$E^{\text{C}}_{r}(\lambda )$ is of the form given in
Theorem~\ref{diffeqthm} with $\omega =\varpi +\varpi^\prime$.

To complete the proof of Theorem~\ref{diffeqthm} it remains to
show that the differential operator $D^{\text{C}}_{r }$
is indeed of the form described by the first part of the theorem
(the commutativity of $D_{1}^{\text{C}},\ldots ,D_{n}^{\text{C}}$
is then an immediate consequence of the commutativity (see
\cite{die:properties}) of
$D_{1,\beta }^{\text{C}},\ldots ,D_{n, \beta}^{\text{C}}$).
More precisely, we should
demonstrate that $D^{\text{C}}_{r }$ is the leading part of
the difference operator $D^{\text{C}}_{r ,\beta }$ in the
formal expansion in $\beta$
and also
that its leading symbol is given by the $r$th elementary symmetric
function in the partials $-\partial^2/\partial x_1^2,\ldots
,-\partial^2/\partial x_n^2$.

To this end we first observe that we may replace
$p_{\lambda ,\beta }^{\text{C}}$ in the l.h.s. of \eqref{limdiffbasis}
by $p_{\lambda }^{\text{C}}$ (using once more that
$\lim_{\beta\rightarrow 0} p_{\lambda ,\beta }^{\text{C}}=p_{\lambda
}^{\text{C}}$).
Moreover, since the polynomials
$\{p_{\lambda }^{\text{C}}\}_{\lambda\in\Lambda}$ constitute
a basis for the space of (even) symmetric polynomials,
we have in fact that
$(D^{\text{C}}_{r,\beta } p)(x)=\beta^{2r} (D^{\text{C}}_{r }p)(x)
+ O(\beta^{2r+1})$ for arbitrary (even) symmetric
polynomial $p(x)$ in the variables $x_1,\ldots ,x_n$.
But then the same holds true for arbitrary
analytic (not necessarily symmetric or even)
function of $x_1,\ldots ,x_n$ and we have (formally)
\begin{equation}\label{diffexp}
D^{\text{C}}_{r,\beta } =\beta^{2r}D^{\text{C}}_{r }
\;\;\; + O(\beta^{2r+1}).
\end{equation}
Here we have used the fundamental property
that the vanishing of a (linear) partial differential operator
on the space of (even) symmetric polynomials implies it be zero
on an arbitrary (analytic) function (i.e., all its coefficients must be zero).
In \cite[Appendix C]{die:commuting} a proof of this general property was
given in a trigonometric context. (Specifically,
we assumed there that the differential operator vanishes
on the space of symmetric
polynomials in $\sin^2(x_1),\ldots ,\sin^2(x_n)$.)
The present case follows after an appropriate substitution of the variables
turning the relevant space of trigonometric polynomials
into the space of (even) symmetric polynomials in $x_1,\ldots ,x_n$.
(Specifically, the type A case is recovered by the substitution
$\sin^2(x_j)\rightarrow x_j$ and the type B case by
$\sin^2(x_j)\rightarrow x_j^2$.)
A consequence of this general property is
that the leading part of the formal
expansion of the difference operator
$D^{\text{C}}_{r,\beta }$ in $\beta$ is completely determined
by its action on the (even) symmetric polynomials. (This excludes the (a
priori)
possibility of a lower-order leading part in the
formal expansion \eqref{diffexp} corresponding to
a term determined by a nontrivial differential
operator that vanishes on the space of (even) symmetric polynomials.)
Notice also that this property
may be used to arrive at an alternative proof for the
commutativity of $D^{\text{C}}_{1 },\ldots ,D^{\text{C}}_{n }$. Indeed,
the commutators $[D^{\text{C}}_{r }, D^{\text{C}}_{s }]$ obviously
vanish on the simultaneous
eigenbasis $\{ p_\lambda^{\text{C}} \}_{\lambda\in \Lambda}$
(and hence on all (even) symmetric polynomials), from which it
then follows
that they be zero identically.

To determine the highest-order symbol of the leading
differential operator $D^{\text{C}}_{r }$ in \eqref{diffexp},
we use that the functions
$v^{\text{C}}$, $w^{\text{C}}$ governing the coefficients of the
difference operator $D^{\text{C}}_{r,\beta }$ are of the
form $1 + O(\beta )$ and that for $v^{\text{C}}, w^{\text{C}}=1$
\begin{eqnarray} \label{lead}
 D_{r,\beta}^{\text{C}}
&\stackrel{(v^{\text{C}},w^{\text{C}}=1)}{=}&
\sum_{\stackrel{J\subset \{ 1,\ldots ,n\}}{|J|= r}}
\prod_{j\in J} \left(
e^{\frac{\beta}{i}\frac{\partial}{\partial x_j} } +
e^{-\frac{\beta}{i}\frac{\partial}{\partial x_j} } -2 \right) \\
&= & \beta^{2r} \Bigl( (-1)^r\!\!\!
\sum_{\stackrel{J\subset \{ 1,\ldots ,n\}}{|J|= r}}
\prod_{j\in J} \frac{\partial^2}{\partial x_j^2} \Bigr) \;\;\;
+ O(\beta^{2r+1}) . \nonumber
\end{eqnarray}
(Notice in this connection that in the situation where
$v^{\text{C}}$ and $w^{\text{C}}$ are taken to be equal to one, the
operators $D_{r,\beta}^{\text{C}}$ \eqref{diffAbetar}, \eqref{diffBbetar}
reduce to
\begin{equation*}
D_{r,\beta}^{\text{A}}
=\sum_{\stackrel{J_+,J_-\subset \{ 1,\ldots ,n\}}
               {J_+\cap J_- =\emptyset ,\; |J_+|+|J_-|\leq r}}
(-2)^{r-|J_+|-|J_-|}
\left( \begin{array}{c} n-|J_+|-|J_-| \\ r-|J_+|-|J_-| \end{array} \right)
e^{\frac{\beta}{i}(\partial_{J_+}-\partial_{J_-})}
\end{equation*}
and
\begin{equation*}
D_{r,\beta}^{\text{B}}=
\sum_{\stackrel{J\subset \{ 1,\ldots ,n\} ,\, 0\leq |J|\leq r}
               {\varepsilon_j=\pm 1,\; j\in J}}
(-2)^{r-|J|} \left( \begin{array}{c} n-|J| \\ r-|J| \end{array} \right)
e^{\frac{\beta}{i}\partial_{\varepsilon J}},
\end{equation*}
which both can be rewritten in the form given by the first line of
\eqref{lead}.)
It then follows (i.e. from the
asymptotics $v^{\text{C}}$, $w^{\text{C}}=1+ O(\beta)$
and \eqref{lead}) that the leading part
$D_{r}^{\text{C}}$ in \eqref{diffexp}
is of the form
\begin{equation*}
D_{r}^{\text{C}}= (-1)^r \sum_{\stackrel{J\subset \{ 1,\ldots ,n\}}{|J|=r}}
\prod_{j\in J}
\frac{\partial^2}{\partial x_j^2}\;\;\; + \text{l.o.}
\end{equation*}
as advertised.

Thus far we have shown that the statements
of Theorem~\ref{diffeqthm} hold when taking
difference operators $D_{r,\beta}^{\text{C}}$ \eqref{diffAbetar},
\eqref{diffBbetar} with $w^{\text{C}}(z)$ given by
\eqref{wArep}, \eqref{wBrep}, where
$\varpi +\varpi^\prime=\omega$ and
$\text{g}_1+\text{g}_1^\prime =g_1$.
The formulation in Theorem~\ref{diffeqthm}
corresponds to choosing the specialization
$\varpi =\omega$, $\varpi^\prime =0$ and
$\text{g}_1=g_1$, $\text{g}_1^\prime =0$.

The symmetry of the differential operator
$D_{r,0}^{\text{C}}$
with respect to the inner product $\langle \cdot ,\cdot \rangle_{\text{C}}$
in the space of permutation-invariant (type $A$) or permutation-invariant
and even (type $B$) polynomials
(Corollary~\ref{symmetrycor})
stems from the fact that the operator in question is diagonal on an orthogonal
basis
(viz. $\{ p_\lambda^{\text{C}} \}_{\lambda\in\Lambda}$) with eigenvalues
that are real.

\subsection{Pieri type recurrence formulas}\label{sec4.3}
In \cite{die:properties} Pieri formulas for the multivariable
continuous Hahn and Wilson polynomials associated
to the weight functions $\Delta^{\text cH}(x)$ \eqref{weightcH}
and $\Delta^{\text W}(x)$ \eqref{weightW} were introduced.
After rescaling and reparametrizing in accordance with
\eqref{scale}, \eqref{repa}, one arrives at the corresponding
Pieri formulas for the polynomials
$p_{\lambda ,\beta}^{\text A}$ and $p_{\lambda ,\beta}^{\text B}$
(which---recall---are determined by the conditions A.1, A.2 and B.1, B.2 of
Section~\ref{sec2.1} with the weight functions
$\Delta^{\text{A}}(x)$ \eqref{weightA} and
$\Delta^{\text{A}}(x)$ \eqref{weightB}
replaced by $\Delta^{\text{A}}_\beta(x)$ \eqref{weightAbeta}
and $\Delta^{\text{B}}_\beta(x)$ \eqref{weightBbeta}).
In the simplest case the resulting Pieri formula takes the form
\begin{eqnarray}\label{recrbeta1}
\hat{E}_{1,\beta}^{\text{C}} (x) P^{\text{C}}_{\lambda ,\beta}(x) &=&
\sum_{\stackrel{1\leq j \leq n}{\lambda +e_j\in\Lambda}}
\hat{V}_{j,\beta }^{\text{C}} (\rho^{\text{C}}+\lambda )
\left( P^{\text{C}}_{\lambda +e_j,\beta}(x) -P^{\text{C}}_{\lambda ,\beta} (x)
\right) +\\
&&
\sum_{\stackrel{1\leq j \leq n}{\lambda -e_j\in\Lambda}}
\hat{V}_{-j,\beta }^{\text{C}} (\rho^{\text{C}}+\lambda )
\left( P^{\text{C}}_{\lambda -e_j,\beta}(x) -P^{\text{C}}_{\lambda ,\beta} (x)
\right)
\nonumber
\end{eqnarray}
with
\begin{equation*}
\hat{V}_{\pm j,\beta }^{\text{C}}(\zeta ) =
\hat{w}^{\text{C}} (\pm \zeta_j)
\prod_{1\leq k\leq n ,k\neq j} \hat{v}^{\text{C}}(\pm \zeta_j +\zeta_k)
\hat{v}^{\text{C}}(\pm \zeta_j -\zeta_k).
\end{equation*}
The functions $\hat{v}^{\text{C}}$, $\hat{w}^{\text{C}}$ are given by
\begin{eqnarray*}
\hat{v}^{\text{A}}(z)&=& 1+\frac{g_0}{z},\;\;\;\;\;
\hat{w}^{\text{A}}(z)=
\frac{(\hat{a}^{\text{A}}+z)(\hat{b}^{\text{A}}+z)}{4z}, \\
\hat{v}^{\text{B}}(z)&=& 1+\frac{g_0}{z},\;\;\;\;\;
\hat{w}^{\text{B}}(z)=
\frac{(\hat{a}^{\text{B}}+z)(\hat{b}^{\text{B}}+z)
(\hat{c}^{\text{B}}+z)(\hat{d}^{\text{B}}+z)}{2z(1+2z)}
\end{eqnarray*}
with
\begin{equation*}
\hat{a}^{\text{A}}= \beta^{-2} (\frac{1}{\varpi}+ \frac{1}{\varpi^\prime})
-1/2,
\;\;\;\;\;\;
\hat{b}^{\text{A}}= \beta^{-2} (\frac{1}{\varpi}- \frac{1}{\varpi^\prime}) +1/2
\end{equation*}
and
\begin{eqnarray*}
&& \hat{a}^{\text{B}}=(\text{g}_1+\text{g}_1^\prime)/2 +
 \beta^{-2} (\frac{1}{\varpi}+ \frac{1}{\varpi^\prime})/2 -1/4,\\
&&\hat{b}^{\text{B}}=(\text{g}_1+\text{g}_1^\prime)/2 -
 \beta^{-2} (\frac{1}{\varpi}+ \frac{1}{\varpi^\prime})/2 +3/4,\\
&& \hat{c}^{\text{B}}=(\text{g}_1-\text{g}_1^\prime)/2 +
 \beta^{-2} (\frac{1}{\varpi}- \frac{1}{\varpi^\prime})/2 +1/4, \\
&&\hat{d}^{\text{B}}=(\text{g}_1-\text{g}_1^\prime)/2 -
 \beta^{-2} (\frac{1}{\varpi}- \frac{1}{\varpi^\prime})/2 +1/4.
\end{eqnarray*}
The vector $\rho^{\text{C}}$ has components given by \eqref{rhoA}, \eqref{rhoB}
and
the multiplyer in the l.h.s. of \eqref{recrbeta1} is given by
\begin{subequations}
\begin{eqnarray}\label{multA}
\hat{E}_{1,\beta}^{\text{A}} (x) &=& -\sum_{1\leq j\leq n} (\frac{ix_j}{\beta}
+\hat{\rho}_j^{\text{A}}) ,\\
\label{multB}
\hat{E}_{1,\beta}^{\text{B}} (x) &=& -\sum_{1\leq j\leq n}
(\frac{x_j^2}{\beta^2}
+(\hat{\rho}_j^{\text{B}})^2)
\end{eqnarray}
\end{subequations}
with
\begin{equation}
\hat{\rho}_j^{\text{A}}= (n-j)g_0 +1/(\beta^2\varpi )
\;\;\;\;\;\;\; \text{and} \;\;\;\;\;\;\;
\hat{\rho}_j^{\text{B}}=(n-j)g_0 +\text{g}_1.
\end{equation}
In the Pieri formula we have furthermore employed the normalization
\begin{equation}\label{Prenobeta}
P_{\lambda ,\beta}^{\text{C}}(x) =c_{\lambda ,\beta}^{\text{C}}\, p_{\lambda
,\beta}^{\text{C}}(x)
\end{equation}
with
\begin{eqnarray*}
c_{\lambda ,\beta}^{\text{A}} &=& (-4i/\beta)^{|\lambda |}
\prod_{1\leq j\leq n} \frac{ [\rho^{\text{A}}_j]_{\lambda_j} }
                           {
[\hat{a}+\rho^{\text{A}}_j,\hat{b}+\rho^{\text{A}}_j]_{\lambda_j}}\\
&& \times
\prod_{1\leq j< k\leq n}
\frac{ [\rho^{\text{A}}_j +\rho^{\text{A}}_k ]_{\lambda_j +\lambda_k} }
     { [g_0+ \rho^{\text{A}}_j +\rho^{\text{A}}_k ]_{\lambda_j +\lambda_k} }
\frac{ [\rho^{\text{A}}_j -\rho^{\text{A}}_k ]_{\lambda_j -\lambda_k} }
     { [g_0+ \rho^{\text{A}}_j -\rho^{\text{A}}_k ]_{\lambda_j -\lambda_k} }
,\\
c_{\lambda ,\beta}^{\text{B}} &=& (-1/\beta^2)^{|\lambda |}
\prod_{1\leq j\leq n} \frac{ [2\rho^{\text{B}}_j]_{2\lambda_j} }
                           {
[\hat{a}+\rho^{\text{B}}_j,\hat{b}+\rho^{\text{B}}_j,
\hat{c}+\rho^{\text{B}}_j,\hat{d}+\rho^{\text{B}}_j ]_{\lambda_j}}\\
&& \times
\prod_{1\leq j< k\leq n}
\frac{ [\rho^{\text{B}}_j +\rho^{\text{B}}_k ]_{\lambda_j +\lambda_k} }
     { [g_0+ \rho^{\text{B}}_j +\rho^{\text{A}}_k ]_{\lambda_j +\lambda_k} }
\frac{ [\rho^{\text{B}}_j -\rho^{\text{B}}_k ]_{\lambda_j -\lambda_k} }
     { [g_0+ \rho^{\text{B}}_j -\rho^{\text{B}}_k ]_{\lambda_j -\lambda_k} } .
\end{eqnarray*}

The recurrence relations of Theorem~\ref{recrthm1} can be
recovered from
\eqref{recrbeta1} for $\beta\rightarrow 0$.
To see this, one first observes that
\begin{equation}\label{PbetatoP}
\lim_{\beta \rightarrow 0} \bigl(\frac{i}{\beta \varpi}\bigr)^{|\lambda |}
P_{\lambda ,\beta}^{\text{A}}(x)=P_{\lambda
}^{\text{A}}(x),\;\;\;\;\;\;\;\;\;\;\;\;
\lim_{\beta \rightarrow 0} P_{\lambda ,\beta}^{\text{B}}(x)=P_{\lambda
}^{\text{B}}(x)
\end{equation}
with $P_{\lambda }^{\text{C}}(x)$ given by \eqref{Preno}.
These limits follow from \eqref{pollim} and the fact that
\begin{equation}\label{cbetatoc}
\lim_{\beta\rightarrow 0} \bigl(\frac{i}{\beta \varpi}\bigr)^{|\lambda |}
c_{\lambda ,\beta}^{\text{A}}= c_{\lambda }^{\text{A}},
\;\;\;\;\;\;\;\;\;\;\;\;
\lim_{\beta\rightarrow 0} c_{\lambda ,\beta}^{\text{B}} = c_{\lambda
}^{\text{B}}
\end{equation}
with $c_{\lambda }^{\text{A}}$ and $c_{\lambda }^{\text{B}}$ given by
\eqref{PrenocA}
and \eqref{PrenocB}. (As usual, we assume an identification of the
parameters of the form
$\omega =\varpi +\varpi^\prime$ and $g_1=\text{g}_1+\text{g}_1^\prime$.)
For type $A$, multiplication of \eqref{recrbeta1} by
$i\beta (\frac{i}{\beta\varpi})^{|\lambda |}$ leads for $\beta\rightarrow 0$
to the first recurrence relation of Theorem~\ref{recrthm1}.
The second (i.e. type $B$) recurrence relation of Theorem~\ref{recrthm1} is
obtained
similarly by sending $\beta$ to zero in
the type $B$ version of \eqref{recrbeta1} after
having multiplied both sides by the factor $(\varpi +\varpi^\prime )\beta^2$.

To derive these limit transitions for the recurrence relations we have used
that for $\beta\rightarrow 0$
\begin{eqnarray*}
&& \hat{v}^{\text{C}}(\rho_j^{\text{C}}+\rho_k^{\text{C}}+\lambda_j+\lambda_k),
\hat{v}^{\text{C}}(-\rho_j^{\text{C}}-\rho_k^{\text{C}}-\lambda_j-\lambda_k)
=  1 + O(\beta^2) ,\\
&&\hat{v}^{\text{C}}(\rho_j^{\text{C}}-\rho_k^{\text{C}}+\lambda_j-\lambda_k)=
1+\frac{g_0}{(k-j)g_0+\lambda_j-\lambda_k}
\end{eqnarray*}
and
\begin{eqnarray*}
&&\hat{w}^{\text{A}}(\rho_j^{\text{A}}+\lambda_j)
= \frac{1}{\beta^2\varpi} ( 1 + O(\beta^2)),\\
&& \hat{w}^{\text{A}}(-\rho_j^{\text{A}}-\lambda_j)
= - \varpi
      \frac{(n-j)g_0+\lambda_j}{2 (\varpi +\varpi^\prime )} ( 1 +
O(\beta^2)),\\
&& \hat{w}^{\text{B}}(\rho_j^{\text{B}}+\lambda_j)=
\frac{(n-j)g_0 +\text{g}_1+\text{g}_1^\prime +1/2+\lambda_j}
{\beta^2 (\varpi +\varpi^\prime)} ( 1 + O(\beta^2)) , \\
&& \hat{w}^{\text{B}}(-\rho_j^{\text{B}}-\lambda_j)=
\frac{(n-j)g_0 +\lambda_j}{\beta^2 (\varpi +\varpi^\prime)}
( 1 + O(\beta^2)).
\end{eqnarray*}
Notice to this end also that in the case of type $A$, a divergent term in the
l.h.s.
of the recurrence relation \eqref{recrbeta1}
originating from the factor $-\sum_j \hat{\rho}_j^{\text{A}}$ (cf.
\eqref{multA})
cancels against
a corresponding divergent term in the r.h.s. originating from the
factor in front
of $P_{\lambda ,\beta}^{\text{A}}$
of the form $-\sum_j \hat{V}_{j,\beta}^{\text{A}}(\rho^{\text{A}}+\lambda)$.
(That the divergent terms on both
sides indeed cancel is seen using the identity
$\sum_j \prod_{k\neq j} ( 1+g_0/(\zeta_j-\zeta_k ))=n$.)

In general the recurrence relations for the polynomials
$P_{\lambda ,\beta}^{\text{C}}$ \eqref{Prenobeta}
induced by \cite{die:properties} become
\begin{eqnarray}\label{recrbetar}
\lefteqn{\hat{E}_{r,\beta}^{\text{C}} (x)\: P_{\lambda ,\beta}^{\text{C}}(x)=}
&&\\
&&\sum_{\stackrel{J\subset \{ 1,\ldots ,n\} ,\, 0\leq|J|\leq r}
               {\varepsilon_j=\pm 1,\; j\in J;\;
                \lambda + e_{\varepsilon J} \in \Lambda}}
\!\!\!\!\!\!\!\!\!
\hat{U}_{J^c,\, r-|J|}^{\text{C}}(\rho^{\text{C}} +\lambda)\,
\hat{V}_{\varepsilon J,\, J^c}^{\text{C}}(\rho^{\text{C}} +\lambda)\,
P_{\lambda +e_{\varepsilon J},\beta}^{\text{C}} (x) ,\nonumber
\end{eqnarray}
$r=1,\ldots ,n$, with
\begin{eqnarray*}
e_{\varepsilon J} &= & \sum_{j\in J} \varepsilon_j e_j
\;\;\;\;\;\;\;\;\;\;\;\;\; (\varepsilon_j \in \{ +1 ,-1\} ),\\
\hat{V}_{\varepsilon J,\, K}^{\text{C}}(\zeta )\!\!\! &=&\!\!\!
\prod_{j\in J} \hat{w}^{\text{C}}(\varepsilon_j\zeta_j)
\prod_{\stackrel{j,j^\prime \in J}{j<j^\prime}}
\hat{v}^{\text{C}}
(\varepsilon_j\zeta_j+\varepsilon_{j^\prime}\zeta_{j^\prime})\,
\hat{v}^{\text{C}}
(\varepsilon_j\zeta_j+\varepsilon_{j^\prime}\zeta_{j^\prime} +1)\\
& & \times
\prod_{\stackrel{j\in J}{k\in K}} \hat{v}^{\text{C}}(\varepsilon_j
\zeta_j+\zeta_k)\,
\hat{v}^{\text{C}}(\varepsilon_j \zeta_j -\zeta_k),\\
\hat{U}_{K,p}^{\text{C}}(\zeta )\!\!\! &=&\\ \!\!\!\!\!\!
 (-1)^p\!\!\!\!\!\!   & &\!\!\!\!\!\!\!\!\!
\sum_{\stackrel{L\subset K,\, |L|=p}
               {\varepsilon_l =\pm 1,\; l\in L }}\!\!
\Bigl( \prod_{l\in L} \hat{w}^{\text{C}}(\varepsilon_l \zeta_l)
\prod_{\stackrel{l,l^\prime \in L}{l<l^\prime}}
\hat{v}^{\text{C}}
(\varepsilon_l\zeta_l+\varepsilon_{l^\prime}\zeta_{l^\prime})\,
\hat{v}^{\text{C}}
(-\varepsilon_l\zeta_l-\varepsilon_{l^\prime}\zeta_{l^\prime} -1 )\\
& &\times
\prod_{\stackrel{l\in L}{k\in K\setminus L}}
\hat{v}^{\text{C}}(\varepsilon_l \zeta_l+\zeta_k)\,
\hat{v}^{\text{C}}(\varepsilon_l \zeta_l -\zeta_k) \Bigr)
\end{eqnarray*}
and
\begin{subequations}
\begin{eqnarray}\label{multipA}
\hat{E}^{\text{A}}_{r,\beta} (x)\!\!\!\! &=&\!\!\!\! (-1)^r
\sum_{\stackrel{ J\subset \{ 1,\ldots ,n\} }{0\leq |J|\leq r}}
 \prod_{j\in J} \frac{ix_j}{\beta}
\sum_{r\leq l_1\leq \cdots \leq l_{r-|J|}\leq n}
\hat{\rho}^{\text{A}}_{l_1}\cdots \hat{\rho}^{\text{A}}_{l_{r-|J|}} ,\\
\hat{E}^{\text{B}}_{r,\beta} (x)\!\!\!\! &=&\!\!\!\! (-1)^{r}
\sum_{\stackrel{ J\subset \{ 1,\ldots ,n\} }{0\leq |J|\leq r}}
 \prod_{j\in J} \frac{x_j^2}{\beta^2}
\sum_{r\leq l_1\leq \cdots \leq l_{r-|J|}\leq n}
(\hat{\rho}^{\text{B}}_{l_1}\cdots \hat{\rho}^{\text{B}}_{l_{r-|J|}})^2 .
\end{eqnarray}
\end{subequations}
For $r=1$ the recurrence formula in \eqref{recrbetar} specializes to that of
\eqref{recrbeta1}.

It is not difficult to see that
the recurrence relations for the multivariable Laguerre polynomials
characterized by
Theorem~\ref{recrthm2} and Theorem~\ref{recrthm3} follow from the
type $B$ version of \eqref{recrbetar} for $\beta\rightarrow 0$.
Indeed, multiplication of \eqref{recrbetar} by the factor $\beta^{2r}(\varpi
+\varpi^\prime )^r$
and sending $\beta$ to zero readily leads to the Laguerre type recurrence
relations.
The verification of this assertion hinges on the second
limit formula of \eqref{PbetatoP}, the
asymptotics
for $\hat{v}^{\text{B}}$, $\hat{w}^{\text{B}}$ displayed above, and
the fact that
\begin{equation*}
\lim_{\beta\rightarrow 0}\beta^{2r}\hat{E}_{r,\beta}^{\text{B}}(x)=
(-1)^r \sum_{\stackrel{ J\subset \{ 1,\ldots ,n\} }{|J|= r}}
 \prod_{j\in J} x_j^2.
\end{equation*}

For type $A$ the transition $\beta\rightarrow 0$ is substantially
more complicated due to the singular nature of the terms in \eqref{recrbetar}.
Specifically,
the multiplyer $\hat{E}_{r,\beta}^{\text{A}}(x)$ \eqref{multipA}
consists of a linear combination
of the elementary symmetric functions in $x_1,\ldots ,x_n$ up to degree $r$.
The coefficients in this linear combination
have a pole at $\beta =0$, the order of which is reversely proportional
to the degree of the elementary symmetric function in question (notice that
$\hat{\rho}_j^{\text{A}}= O(\beta^{-2})$).
Hence, for $\beta\rightarrow 0$ the contributions of the lower-degree
elementary symmetric functions to $\hat{E}_{r,\beta}^{\text{A}}(x)$
become predominant.
To get rid of these lower-degree divergent terms, we
take an appropriate linear combination of
the recurrence relations \eqref{recrbetar}
that cast them into a system of the form
\begin{equation}\label{recrrew}
\Bigl( \sum_{\stackrel{J\in \{ 1,\ldots ,n\} }{|J|=r} } \prod_{j\in J} x_j
\Bigr)
\tilde{P}_{\lambda ,\beta}^{\text{A}} (x)=
\sum_{\stackrel{J\subset \{ 1,\ldots ,n\} ,\, 0\leq|J|\leq r}
               {\varepsilon_j=\pm 1,\; j\in J;\;
              e_{\varepsilon J} +\lambda \in \Lambda}}
\!\!\!\!\!\!\!\!\!
\hat{W}_{\varepsilon J,\, J^c;r,\beta}^{\text{A}}\,
\tilde{P}_{\lambda +e_{\varepsilon J},\beta}^{\text{A}} (x)
\end{equation}
with
\begin{equation*}
\tilde{P}_{\lambda ,\beta}^{\text{A}} (x) =
\Bigl( \frac{i}{\beta\varpi} \Bigr)^{|\lambda |} P_{\lambda ,\beta}^{\text{A}}
(x) .
\end{equation*}
Basically, this boils down
to passing from Pieri formulas corresponding to the symmetric functions
$\hat{E}_{r,\beta}^{\text{A}}(x)$ \eqref{multipA} to
Pieri formulas corresponding to the elementary symmetric functions
$\sum_{|J|=r} \prod_{j\in J} x_j$
by subtracting from the $r$th Pieri formula in
\eqref{recrbetar}
a suitable linear combination of the Pieri formulas corresponding
to $\hat{E}_{s,\beta}^{\text{A}}(x)$  with $s<r$ (and multiplication
by an overall factor).
The coefficients of the terms in the r.h.s. of \eqref{recrbetar}
labeled by index sets $J$ with $|J|=r$ are
invariant with respect to such changes in the l.h.s. (up to an
overall factor
$(i\beta)^r$ and a factor caused by the
change of the normalization
$P_{\lambda ,\beta}^{\text{A}}\rightarrow \tilde{P}_{\lambda
,\beta}^{\text{A}}$).
More precisely, we obtain that
for $|J|=r$ the coefficient $\hat{W}_{\varepsilon J,\, J^c;r,\beta}^{\text{A}}$
in the r.h.s. of \eqref{recrrew} is given by
\begin{equation*}
\hat{W}_{\varepsilon J,\, J^c;r,\beta}^{\text{A}} =
 (i\beta)^r (-i\beta \varpi )^{\sum_{j\in J} \varepsilon_j}
\hat{V}_{\varepsilon J,\, J^c}^{\text{A}}(\rho^{\text{A}} +\lambda) .
\end{equation*}
The type $A$ version of Theorem~\ref{recrthm2} then follows for $\beta
\rightarrow 0$
(using the limit formula \eqref{PbetatoP} and the above asymptotics for
$\hat{v}^{\text{A}}$ and $\hat{w}^{\text{A}}$).

It remains to verify the normalization properties of $P_\lambda^{\text{C}}(x)$
stated in Theorem~\ref{normalizationthm}.
These properties are a consequence of the fact that
$P_{\lambda ,\beta}^{\text{C}} (i\beta \hat{\rho}^{\text{C}})=1$ (see
the remarks in \cite[Sec 6]{die:properties}).
Specifically, by sending $\beta$ to zero in the relation
$P_{\lambda ,\beta}^{\text{C}} (i\beta \hat{\rho}^{\text{C}})=1$
we arrive at the normalization properties of Theorem~\ref{normalizationthm}.
For type $B$ this is immediate from the limit formula \eqref{PbetatoP};
for type $A$ this is seen by noticing that
$\lim_{\beta\rightarrow 0} P_{\lambda ,\beta}^{\text{A}} (i\beta
\hat{\rho}^{\text{A}})$
picks up the highest-degree  homogeneous part of $P_{\lambda
,\beta}^{\text{A}}(x)$
evaluated in $x=\mathbf{1}$
(here we use that $\hat{\rho}^{\text{A}}_j = \frac{1}{\varpi \beta^2} + O(1)$
together with the limit formula \eqref{cbetatoc}), which is equal to
$\lim_{\alpha\rightarrow\infty} \alpha^{-|\lambda |}
P_\lambda^{\text{A}}(\alpha \mathbf{1})$.

\setcounter{equation}{0}
\renewcommand{\theequation}{A.\arabic{equation}}
\section*{Appendix: Convergence of the weight functions}
In this appendix it will be shown that the weight functions
$\Delta^{\text{A}}_\beta (x)$ \eqref{weightAbeta} and
$\Delta^{\text{B}}_\beta (x)$ \eqref{weightBbeta}
converge for $\beta\rightarrow 0$ to the weight functions
$\Delta^{\text{A}}(x)$ \eqref{weightA} and
$\Delta^{\text{B}}(x)$ \eqref{weightB}, respectively.
More precisely, we will prove the somewhat stronger result that
\begin{eqnarray}\label{limit}
& &\lim_{\beta\rightarrow 0}\;
\int_{-\infty}^{\infty} \cdots \int_{-\infty}^{\infty}
p(x)\: \Delta^{\text{C}}_\beta (x)\: dx_1\cdots dx_n
 =  \\
&& \makebox[3em]{}
\int_{-\infty}^{\infty} \cdots \int_{-\infty}^{\infty}
p(x)\: \Delta^{\text{C}}(x)\: dx_1\cdots dx_n
\;\;\;\;\;\; (C= A\; \text{or}\; B),  \nonumber
\end{eqnarray}
where $p(x)$ denotes an arbitrary polynomial in the variables $x_1,\ldots ,
x_n$.
In \eqref{limit} it is understood that
the parameters of $\Delta^{\text{C}}(x)$ \eqref{weightA}, \eqref{weightB}
and $\Delta^{\text{C}}_\beta (x)$ \eqref{weightAbeta}, \eqref{weightBbeta}
are related via the identification
$\omega \equiv \varpi +\varpi^\prime$ and $g_1\equiv
\text{g}_1+\text{g}_1^\prime$.

Let us start by inferring the pointwise convergence of the weight functions.
For this purpose we use the limit formula
\begin{equation}\label{lim1}
\lim_{\beta \rightarrow 0}\;
\delta (\alpha ,\beta)\:
|\Gamma ( \frac{1}{\alpha \beta^2 }+ i\beta^{-1} y ) |
 = \exp (-\alpha y^2 /2)\;\;\;\;\;\;\;\;\;\;\;\;\; (\alpha >0)
\end{equation}
with
\begin{equation*}
\delta (\alpha ,\beta )\equiv \sqrt{\frac{e}{2\pi}}\:
e^{(1+\log (\alpha\beta^2)) (\alpha^{-1}\beta^{-2}-1/2)}
\end{equation*}
and the limit formula
\begin{equation}\label{lim2}
 \lim_{\beta \rightarrow 0}\;
\left| \beta^a \frac{\Gamma (a + b +i\beta^{-1}y )}{\Gamma (b + i\beta^{-1}y )}
\right| =
|y|^a  \;\;\;\;\;\;\;\;\;\;\;\;\;\;\;\;\; (a,b \geq 0),
\end{equation}
where in both formulas it is assumed that $y$ and $\beta$ are real.
By applying \eqref{lim1} and \eqref{lim2} to the factors of
$\Delta^{\text{A}}_\beta (x)$ \eqref{weightAbeta} and
$\Delta^{\text{B}}_\beta (x)$ \eqref{weightBbeta}, one readily sees
that for $\beta\rightarrow 0$ these weight functions converge pointwise to
$\Delta^{\text{A}} (x)$ \eqref{weightA}
and $\Delta^{\text{B}} (x)$ \eqref{weightB} as indicated.
The normalization factors of the form $\delta (\alpha ,\beta )$ and
$|\beta |^a$ in \eqref{lim1} and \eqref{lim2}
ensure a finite and nontrivial limit;
the factors in question have been collected in the weight functions
$\Delta^{\text{A}}_\beta (x)$ \eqref{weightAbeta} and
$\Delta^{\text{B}}_\beta (x)$ \eqref{weightBbeta} into
the overall normalization constants
$\mathcal{D}^{\text{A}}(\beta )=|\beta |^{n(n-1)g_0} \delta (\varpi
,\beta)^{2n}
\delta (\varpi^\prime ,\beta )^{2n}$
and $\mathcal{D}^{\text{B}}(\beta )=|\beta
|^{2n(n-1)g_0+2n(\text{g}_1+\text{g}_1^\prime)}
 \delta (\varpi ,\beta)^{2n}$ $\delta (\varpi^\prime ,\beta )^{2n}$.

The limit formulas \eqref{lim1} and \eqref{lim2}
may be verified with the aid of
Stirling's formula for the asymptotics of the gamma function
for large values of the argument, which
reads (see e.g. \cite{abr-ste:handbook,olv:asymptotics})
\begin{equation}\label{stirling}
\Gamma (z) = (2\pi)^{1/2}\, e^{-z} z^{z-1/2}\cdot \exp(R(z))
\end{equation}
with $R(z) = O(1/|z|)$ for $|z|\rightarrow \infty$ in the sector $|\arg (z)| <
\pi$.
Substitution of $z = \alpha^{-1}\beta^{-2} +i\beta^{-1}y$ in \eqref{stirling}
entails that for $\beta\rightarrow 0$
\begin{eqnarray}\label{as1}
\lefteqn{|\Gamma ( \frac{1}{\alpha \beta^2 }+ i\beta^{-1} y ) | =} &&\\
&& \frac{1}{\delta (\alpha ,\beta )}
(1+\alpha^2\beta^2y^2)^{\frac{1}{2\alpha\beta^2}}\, e^{-\frac{y}{\beta}\arctan
(\alpha \beta y )}
\, (1+O(\beta^2 ) ), \nonumber
\end{eqnarray}
which implies \eqref{lim1}.
In a similar way one concludes from \eqref{stirling} that for $\beta\rightarrow
0$
\begin{equation}\label{as2}
\left| \frac{\Gamma (a + b +i\beta^{-1}y )}{\Gamma (b + i\beta^{-1}y )} \right|
=
| y/\beta |^a \, (1+ O(\beta ))  ,
\end{equation}
which implies \eqref{lim2}.

Let us next demonstrate that the pointwise convergence of
the integrands carries over to the convergence of the integrals
by invoking Lebesgue's dominated convergence theorem.
For this purpose it is needed to dominate
(the absolute value of) the integrand in the l.h.s. of \eqref{limit} uniformly
in $\beta$ by an integrable function.
We will do so by bounding individually the factors comprising the weight
function.
Specifically, it turns out that factors of the form
$\delta (\varpi ,\beta)\:
|\Gamma ( \frac{1}{\varpi \beta^2 }+ i\beta^{-1} x_j ) |$
(and $\delta (\varpi^\prime ,\beta)\:
|\Gamma ( \frac{1}{\varpi^\prime \beta^2 }+ i\beta^{-1} x_j ) |$ )
may be dominated by an exponentially decaying function
and that the remaining factors---which consist of ratios
of the form $|\beta^{g_0}\Gamma (g_0 +i\beta^{-1}(x_j \pm x_k))/\Gamma
(i\beta^{-1}(x_j \pm x_k))|$,
$|\beta^{\text{g}_1}\Gamma (\text{g}_1 +i\beta^{-1}x_j )/\Gamma (i\beta^{-1}x_j
)|$ and
$|\beta^{\text{g}_1^\prime }\Gamma (\text{g}_1^\prime +1/2 +i\beta^{-1}x_j )/
\Gamma (1/2+i\beta^{-1}x_j )|$---grow at most polynomially in the variables
$x_1,\ldots ,x_n$.
Hence, the integrand on the l.h.s. may be dominated by an exponentially
decaying function and \eqref{limit} follows from the pointwise convergence
of the integrands by the dominated convergence theorem as advertised.

To validate the above claims regarding the bounds on the growth of the factors
constituting
the weight functions $\Delta^{\text{A}}_\beta (x)$ \eqref{weightAbeta} and
$\Delta^{\text{B}}_\beta (x)$ \eqref{weightBbeta},
we need a precise estimate for the error term $R(z)$
appearing in the Stirling formula \cite{olv:asymptotics}
\begin{equation}\label{errorbound}
|R(z)| \leq \frac{1}{12 |z| \cos^2 (\theta/2 ) }\;\;\;\;\;\;\;\;\; (\theta =
\arg (z)).
\end{equation}
(This estimate is valid in the whole sector $|\arg (z)| < \pi$.)
By substituting $z = \alpha^{-1}\beta^{-2} +i\beta^{-1}y$
in the Stirling formula \eqref{stirling}
(where $\alpha $ is assumed to be positive) and
combining with the error estimate \eqref{errorbound}, we find (cf. \eqref{as1})
\begin{equation}\label{bound1}
\delta (\alpha ,\beta)\:
|\Gamma ( \frac{1}{\alpha \beta^2 }+ i\beta^{-1} y ) | \leq
e^{-F_\beta (y)}G_\beta ,
\end{equation}
with
\begin{equation}
F_\beta (y) = \frac{y}{\beta}\arctan (\alpha \beta y )-
              \frac{1}{2\alpha \beta^2}\log (1+\alpha^2 \beta^2y^2)
\end{equation}
and $G_\beta = e^{\alpha \beta^2 /6}$.
(In our situation $\cos^2 (\theta /2) \geq 1/2$ in view of the fact that
the real part of $z=\frac{1}{\alpha \beta^2 }+ i\beta^{-1} y$
is positive.) Differentiation of $F_\beta (y)$ with respect to $y$ yields
\begin{equation}
\partial_y  F_\beta (y)= \beta^{-1}\arctan (\beta\alpha y),
\end{equation}
which shows that $F_\beta (y)$ is nonnegative as
an increasing/decreasing function for $y$ positive/negative with $F_\beta
(0)=0$.
{}From the asymptotics for
$|y|\rightarrow \infty$ one furthermore sees that
the factor $\exp (-F_\beta(y))$ decays exponentially.
A little more precise analysis reveals that for $0 < \beta < 1$
\begin{equation}\label{expobound}
e^{-F_\beta (y)} \leq \left\{ \begin{array}{ll}
                        1 & \text{for}\; |y|< 1/\alpha \\
                        e^{-|y|/3} & \text{for}\; |y| \geq 1/\alpha .
                        \end{array} \right.
\end{equation}
To obtain the exponential bound on the tail
we have used: (i.) that for $0<\beta <1$
\begin{eqnarray*}
F_\beta (1/\alpha ) &=& \frac{1}{\alpha\beta}\arctan ( \beta  )-
              \frac{1}{2 \alpha \beta^2}\log (1+ \beta^2) \\
              & > & ( \pi/4 - \log(\sqrt{2}) )/\alpha > 1/(3 \alpha ),
\end{eqnarray*}
(ii.) that for $y\geq 1/\alpha$ and $0<\beta < 1$
\begin{equation*}
\left(\partial_y  F_\beta \right) (y) \geq
\left(\partial_y  F_\beta \right) (1/\alpha ) = \beta^{-1}\arctan (\beta) > \pi
/4 > 1/3,
\end{equation*}
and (iii.) that $F_\beta (y)$ is even in $y$.
We conclude from \eqref{bound1} and \eqref{expobound} that for $0<\beta <1$
the factors in the weight function $\Delta_\beta^{\text{C}}(x)$ ($C=A,B$) of
the form
$\delta (\varpi ,\beta)\:
|\Gamma ( \frac{1}{\varpi \beta^2 }+ i\beta^{-1} x_j ) |$
(and $\delta  (\varpi^\prime ,\beta)\:
|\Gamma ( \frac{1}{\varpi^\prime \beta^2 }+ i\beta^{-1} x_j ) |$ )
may indeed be dominated uniformly in $\beta$ by an exponentially decaying
function of $x_j$.

It remains to check that the rest of the
factors---which consist of ratios of gamma functions---can be dominated
by a function that grows at most polynomially in the variables.
To this end we should find bounds on the gamma
function ratios of the type appearing in the l.h.s. of \eqref{lim2}.
Notice that for $a\in \mathbb{N}$ we have
\begin{equation}
\left| \beta^a \frac{\Gamma (a + b +i\beta^{-1}y )}{\Gamma (b + i\beta^{-1}y )}
\right|
= \prod_{m=0}^{a-1} |(m+b)\beta +iy |,
\end{equation}
which is easily dominated uniformly in $0<\beta <1$ by a function with
polynomial growth in $y$ (take
e.g. the function $((a+b)^2 +y^2)^{a/2}$).
The case of general positive (not necessarily integer valued) $a$ is a little
less straightforward;
it will be addressed here with the aid of the integral representation
\begin{equation}\label{ratio}
\frac{\Gamma (a+z)}{\Gamma (z)} = z^a
\exp \left( \int_0^\infty  e^{-zt} \left ( a -\frac{1-e^{-at}}{1-e^{-t}}\right)
\frac{dt}{t} \right)
\;\;\;\;\;\;\;\;\;\;\;\;\;\; \text{Re}(z) > 0,
\end{equation}
which can be obtained by integrating Gauss' integral formula for the psi
function
$\psi (z) = \Gamma^\prime (z) /\Gamma (z)$
\cite{abr-ste:handbook,olv:asymptotics}
\begin{equation}
\psi (z) = \int_0^\infty  \left ( \frac{e^{-t}}{t}
-\frac{e^{-zt}}{1-e^{-t}}\right) dt
\;\;\;\;\;\;\;\;\;\;\;\;\;\; \text{Re}(z) > 0
\end{equation}
(using also that
$\log (z)=\int_0^\infty t^{-1}( e^{-t}-e^{-tz} ) dt$
for $\text{Re}(z) > 0$).
Let us for the moment assume that $b$ is positive. Then substitution of $z = b
+i\beta^{-1}y$
in \eqref{ratio} entails
\begin{eqnarray}
\lefteqn{ \left| \beta^a  \frac{\Gamma (a+b+i\beta^{-1}y )}{\Gamma
(b+i\beta^{-1}y )} \right| =} &&\\
&& (\beta^2 b^2 + y^2 )^{a/2}
\exp \left( \int_0^\infty  e^{-bt}\cos ( yt/\beta )
 \left ( a -\frac{1-e^{-at}}{1-e^{-t}}\right) \frac{dt}{t} \right) .\nonumber
\end{eqnarray}
The integral within the exponent is bounded by a constant
with value $$\int_0^\infty  e^{-bt} \left | a
-\frac{1-e^{-at}}{1-e^{-t}}\right| \frac{dt}{t}$$
and the factor in front is smaller than $(b^2 +y^2)^{a/2}$ for $0<\beta <1$.
The case that $b$ becomes zero can be reduced to the previous situation with
positive $b$ by
means of the identity
\begin{equation}
\frac{\Gamma (a +i\beta^{-1}y)}{\Gamma (i\beta^{-1}y)} =
 \frac{\Gamma (a +1+i\beta^{-1}y)}{\Gamma (1+i\beta^{-1}y)}
\left(\frac{iy}{a\beta + iy}\right) .
\end{equation}
The upshot is that
for $0<\beta <1$ the factors
$|\beta^{g_0}\Gamma (g_0 +i\beta^{-1}(x_j \pm x_k))/\Gamma (i\beta^{-1}(x_j \pm
x_k))|$,
$|\beta^{\text{g}_1}\Gamma (\text{g}_1 +i\beta^{-1}x_j )/\Gamma (i\beta^{-1}x_j
)|$ and
$|\beta^{\text{g}_1^\prime }\Gamma (\text{g}_1^\prime +1/2 +i\beta^{-1}x_j )/
\Gamma (1/2+i\beta^{-1}x_j )|$
can be uniformly dominated in $\beta $ by a function that grows at most
polynomially
in the variables $x_1,\dots ,x_n$, which completes the proof of \eqref{limit}.

\vfill
\section*{Acknowledgments}
The author would like to thank S.N.M. Ruijsenaars for
some helpful remarks in connection with the material in the appendix
(in particular with respect to
the usefulness of the integral representation \eqref{ratio}).

\bibliographystyle{amsplain}

\begin{thebibliography}{000000}

\bibitem[A1]{ask:some} R. Askey, {\em Some basic hypergeometric extensions
of integrals of Selberg and Andrews}, SIAM J. Math. Anal. {\bf 11} (1980),
938--951.

\bibitem[A2]{ask:continuous} \bysame, {\em Continuous Hahn polynomials},
J. Phys. A: Math. Gen. {\bf 18} (1985), L1017--L1019.

\bibitem[AtSu]{ata-sus:hahn} N. M. Atakishiyev and S. K. Suslov,
{\em The Hahn and Meixner polynomials of an imaginary argument and
some of their applications}, J. Phys. A: Math. Gen. {\bf 18} (1985),
1583--1596.

\bibitem[AbSt]{abr-ste:handbook} M. Abramowitz and I. A. Stegun (eds.), {\em
Handbook
of mathematical functions}, Dover Publications, New York, 1972 (9th printing).

\bibitem[BF]{bak-for:calogero-sutherland} T. H. Baker and P. J. Forrester,
{\em The Calogero-Sutherland model and generalized classical polynomials},
Preprint Research Institute for Mathematical
Sciences, Kyoto, RIMS-1094, 1996.

\bibitem[BO]{bee-opd:certain} R. J. Beerends and E. M. Opdam, {\em Certain
hypergeometric series
related to the root system $BC$}, Trans. Amer. Math. Soc. {\bf 339} (1993),
581--609.

\bibitem[BHV]{bri-han-vas:explicit} L. Brink, T. H. Hansson, and M. A.
Vasiliev,
{\em Explicit solution of the n-body Calogero model},
Phys. Lett. B {\bf 286} (1992), 109--111.

\bibitem[BHKV]{bri-han-kon-vas:calogero} L. Brink, T. H. Hansson, S. Konstein,
and
M. A. Vasiliev, {\em The Calogero model---anyonic representation, fermionic
extension and supersymmetry}, Nucl. Phys. B {\bf 401} (1993), 591--612.

\bibitem[Ca]{cal:solution} F. Calogero, {\em Solution of the one-dimensional
$n$-body problems
with quadratic and/or inversely quadratic pair potentials},
J. Math. Phys. {\bf 12} (1971), 419--436.

\bibitem[Co]{con:distribution} A. G. Constantine, {\em The distribution of
Hotelling's
generalized $T_0^2$}, Ann. Math. Statist. {\bf 37} (1966), 215--225.


\bibitem[De1]{deb:polynomes} A. Debiard,
{\em Polyn\^omes de Tch\'ebychev et de Jacobi dans un espace euclidien
de dimension $p$},
C. R. Acad. Sci. Paris S\'er. I Math. {\bf 296} (1983), 529--532.

\bibitem[De2]{deb:systeme} \bysame, {\em Syst\`eme diff\'erentiel
hyperg\'eom\'etrique
et parties radiales des op\'erateurs invariants des espaces sym\'etriques de
type $BC_p$},
in: S\'eminaire d'Alg\`ebre Paul Dubreil et Marie-Paule Malliavin (M.-P.
Malliavin, ed.),
Lecture Notes in Math., vol. 1296, Springer, Berlin, 1988, pp. 42--124.

\bibitem[D1]{die:commuting} J. F. van Diejen, {\em Commuting difference
operators
with polynomial eigenfunctions}, Compositio Math. {\bf 95} (1995), 183--233.

\bibitem[D2]{die:difference} \bysame, {\em Difference Calogero-Moser systems
and finite Toda chains}, J. Math. Phys. {\bf 36} (1995), 1299--1323.

\bibitem[D3]{die:multivariable} \bysame, {\em Multivariable continuous Hahn
and Wilson polynomials related to integrable difference systems},
J. Phys. A: Math. Gen. {\bf 28} (1995), L369--L374.

\bibitem[D4]{die:properties} \bysame, {\em Properties of some
families of hypergeometric orthogonal polynomials in several
variables}, Math. preprint University of Tokyo UTMS 96-10, 1996.

\bibitem[Du1]{dun:orthogonal} C. F. Dunkl, {\em Orthogonal polynomials on the
sphere
with octahedral symmetry}, Trans. Amer. Math. Soc. {\bf 282} (1984), 555--575.

\bibitem[Du2]{dun:reflection} \bysame, {\em Reflection groups and orthogonal
polynomials on the sphere}, Math. Z. {\bf 197} (1988), 33--60.

\bibitem[Du3]{dun:differential-difference} \bysame, {\em
Differential-difference
operators associated to reflection groups}, Trans. Amer. Math. Soc. {\bf 311}
(1989), 167--183.

\bibitem[G]{gam:exact} P. J. Gambardella, {\em Exact results in quantum
many-body
systems of interacting particles in many dimensions with
$\overline{SU(1,1)}$ as the dynamical group}, J. Math. Phys. {\bf 16} (1975),
1172--1187.

\bibitem[H]{hec:elementary} G. J. Heckman, {\em An elementary approach to the
hypergeometric shift operator of Opdam}, Invent. Math. {\bf 103} (1991),
341--350.

\bibitem[He]{her:bessel} C. S. Herz, {\em Bessel functions of matrix argument},
Ann. Math. {\bf 61} (1955), 474--523.

\bibitem[J]{jam:special} A. T. James, {\em Special functions of matrix and
single argument in statistics}, in: Theory and applications of special
functions (R. Askey, ed.), Academic Press, New York, 1975, pp. 497--520.

\bibitem[Ka]{kak:common} S. Kakei, {\em Common algebraic structure
for the Calogero-Sutherland models}, Preprint 1996.

\bibitem[KS]{koe-swa:askey-scheme} R. Koekoek and R. F. Swarttouw, {\em The
Askey-scheme of hypergeometric orthogonal polynomials and its q-analogue},
Math. report Delft University of Technology 94-05, 1994.

\bibitem[LV1]{lap-vin:rodrigues} L. Lapointe and L. Vinet, {\em A Rodrigues
formula for the Jack polynomials and the Macdonald-Stanley conjecture},
Internat. Math. Res. Notices {\bf 1995}, 419--424.

\bibitem[LV2]{lap-vin:exact} \bysame, {\em Exact operator solution of the
Calogero-Sutherland model}, Commun. Math. Phys. {\bf 178} (1996), 425--455.

\bibitem[La1]{las:formule} M. Lassalle, {\em Une formule du bin\^ome
g\'en\'eralis\'ee
pour les polyn\^omes de Jack},
C. R. Acad. Sci. Paris S\'er. I Math. {\bf 310} (1990), 253--256.

\bibitem[La2]{las:laguerre} \bysame, {\em Polyn\^omes de Laguerre
g\'en\'eralis\'es},
C. R. Acad. Sci. Paris S\'er. I Math. {\bf 312} (1991), 725--728.

\bibitem[La3]{las:hermite} \bysame, {\em Polyn\^omes de Hermite
g\'en\'eralis\'es},
C. R. Acad. Sci. Paris S\'er. I Math. {\bf 313} (1991), 579--582.

\bibitem[M1]{mac:some} I. G. Macdonald, {\em Some conjectures for root
systems},
SIAM J. Math. Anal. {\bf 13} (1982), 988--1007.

\bibitem[M2]{mac:commuting} \bysame, {\em Commuting differential operators
and zonal spherical functions}, in: Algebraic groups Utrecht 1986
(A. M. Cohen, W. H. Hesselink, W. L. J. van der Kallen,
and J. R. Strooker, eds.), Lecture Notes in Math., vol. 1271, Springer, Berlin,
1987, pp. 189--200.

\bibitem[M3]{mac:hypergeometric} \bysame, {\em Hypergeometric functions},
unpublished manuscript.

\bibitem[M4]{mac:symmetric} \bysame, {\em Symmetric functions and Hall
polynomials}, 2nd edition, Oxford mathematical monographs,
Clarendon Press, Oxford, 1995.

\bibitem[Me]{meh:random} M. L. Mehta, {\em Random matrices}, 2nd edition,
Academic Press, Boston, 1991.

\bibitem[Mu]{mui:aspects} R. J. Muirhead, {\em Aspects of multivariate
statistical theory}, Wiley, New York, 1982.

\bibitem[OOS]{och-osh-sek:commuting} H. Ochiai, T. Oshima, and H. Sekiguchi,
{\em Commuting families of symmetric differential operators},
Proc. Japan Acad. Ser. A Math. Sci. {\bf 70} (1994), 62--66.

\bibitem[OP]{ols-per:quantum} M. A. Olshanetsky and A. M. Perelomov,
{\em Quantum integrable systems related to Lie algebras},
Phys. Rep. {\bf 94} (1983), 313--404.

\bibitem[Ol]{olv:asymptotics} F. W. J. Olver, {\em Asymptotics and
special functions}, Academic Press, New York, 1974.

\bibitem[OS]{osh-sek:commuting} T. Oshima and H. Sekiguchi, {\em Commuting
families of differential operators invariant under the action of a
Weyl group}, J. Math. Sci. Univ. Tokyo {\bf 2} (1995), 1--75.

\bibitem[Pe]{per:algebraic} A. M. Perelomov, {\em Algebraic approach to the
solution
of a one-dimensional model of $n$ interacting particles},
Theoret. and Math. Phys. {\bf 6} (1971), 263--282.

\bibitem[Po]{pol:exchange} A. P. Polychronakos, {\em Exchange operator
formalism for integrable systems of particles}, Phys. Rev. Lett.
{\bf 69} (1992), 703--705.

\bibitem[R]{rui:complete} S. N. M. Ruijsenaars, {\em Complete integrability of
relativistic Calogero-Moser systems and elliptic function identities},
Commun. Math. Phys. {\bf 110} (1987), 191--213.

\bibitem[S]{sek:zonal} J. Sekiguchi, {\em Zonal spherical functions on some
symmetric spaces}, Publ. RIMS Kyoto Univ. {\bf 12} Suppl. (1977), 455--459.

\bibitem[Se]{sel:bemerkninger} A. Selberg, {\em Bemerkninger om et multipelt
integral},
Norsk Mat. Tidsskr. {\bf 26} (1944), 71--78
(Collected papers, vol. 1, Springer, Berlin, 1989, pp. 204--213).

\bibitem[St]{sta:some} R. P. Stanley, {\em Some combinatorial properties of
Jack symmetric functions}, Adv. Math. {\bf 77} (1989), 76--115.

\bibitem[UW]{uji-wad:rodrigues} H. Ujino and M. Wadati, {\em Rodrigues formula
for hi-Jack symmetric polynomials associated with the quantum Calogero model},
J. Phys. Soc. Japan {\bf 65} (1996), 2423--2439.

\bibitem[V]{vre:formulas} L. Vretare, {\em Formulas for elementary spherical
functions and generalized Jacobi polynomials}, SIAM J. Math. Anal. {\bf 15}
(1984),
805--833.

\bibitem[W]{wil:some} J. A. Wilson, {\em Some hypergeometric orthogonal
polynomials},
SIAM J. Math. Anal. {\bf 11} (1980), 690--701.

\end{thebibliography}

\end{document}